\documentclass{article}
\usepackage[utf8]{inputenc}

\usepackage{geometry}
\usepackage{mathrsfs}
\usepackage{bm} 
\usepackage{bbm}
\usepackage{authblk}
\usepackage{xcolor}
\usepackage{comment}
\usepackage{upgreek}
\usepackage{amssymb}  
\usepackage{amsthm}
\usepackage{amsmath}
\usepackage{amsfonts}
\usepackage{graphicx}
\usepackage{pdfpages}

\newtheorem*{remark}{Remark}
\usepackage{algorithm}
\usepackage{algorithmic}
\usepackage{tabularx}

\usepackage{color}
\definecolor{deepblue}{rgb}{0,0,0.5}
\definecolor{deepred}{rgb}{0.6,0,0}
\definecolor{deepgreen}{rgb}{0,0.5,0}

\usepackage{listings}

\providecommand{\keywords}[1]
{
	\textbf{\textit{Keywords}~} #1
}

\definecolor{codegreen}{rgb}{0,0.6,0}
\definecolor{codegray}{rgb}{0.5,0.5,0.5}
\definecolor{codepurple}{rgb}{0.58,0,0.82}
\definecolor{backcolour}{rgb}{0.95,0.95,0.92}

\lstdefinestyle{mystyle}{
    backgroundcolor=\color{backcolour},   
    commentstyle=\color{codegreen},
    keywordstyle=\color{magenta},
    numberstyle=\tiny\color{codegray},
    stringstyle=\color{codepurple},
    basicstyle=\ttfamily\small,
    breakatwhitespace=false,         
    breaklines=true,                 
    captionpos=b,                    
    keepspaces=true,                 
    numbers=none,                    
    numbersep=5pt,                  
    showspaces=false,                
    showstringspaces=false,
    showtabs=false,                  
    tabsize=2
}
 
\usepackage{lmodern}
\usepackage{upquote}
\usepackage{fancyvrb}

\newcommand\textcode{\Verb}

\def\correspondingauthor{\footnote{Corresponding author.}}

\graphicspath{{figures/}}

\title{A finite element method for simulating soft active non-shearable rods immersed in generalized Newtonian fluids}
\author{Roberto Federico Ausas$^1$, Cristian Guillermo Gebhardt$^{2,}$\correspondingauthor{} , Gustavo Carlos Buscaglia$^1$}

\date{$^1$ Instituto de Ciências Matemáticas e de Computação, Universidade de São Paulo, Av. Trabalhador São-carlense, 400, 13566-590, São Carlos, SP, Brazil\\
		~\\
		$^2$ Geophysical Institute (GFI) and Bergen Offshore Wind Centre (BOW), University of Bergen, Allégaten 70, 5007 Bergen, Norway\\
	     ~\\
		\today}

\newcommand{\field}{\mathbb}

\newcommand{\reals}{\field{R}}

\newcommand{\ie}{\textit{i.e.,}~}

\newcommand{\thematrix}[1]{\mathbf{#1}}
\newcommand{\thevector}[1]{\bm{#1}}
\newcommand{\theoperator}[1]{\mathscr{#1}}

\newcommand{\posname}{m}
\newcommand{\pos}{\mathbf{\posname}}

\newcommand{\dspos}{\pos'}

\newcommand{\dsdspos}{\pos''}
\newcommand{\varpos}{\delta\pos}

\newcommand{\eye}{\thematrix{I}}

\newcommand{\Pos}{\thevector{\varPhi}}

\newcommand{\dxiPos}{\pd{\xi}\Pos}
\newcommand{\dxidxiPos}{\pd{\xi}^2\Pos}

\newcommand{\moddspos}{|\dspos|}

\newcommand{\moddxiPos}{|\dxiPos|}

\newcommand{\pd}[1]{\partial_{#1}}
\newcommand{\pairing}{\cdot}

\newcommand{\jye}{\thematrix{J}}

\newcommand{\tanname}{t}

\newcommand{\tanv}{\mathbf{\check{t}}}%{\thevector{\tanname}}
\newcommand{\norv}{\mathbf{\check{n}}}%{\thevector{\norname}}

\newcommand{\Projnortanv}{\mathfrak{P}^{\perp}_{\tanv}}

\newcommand{\interface}{\Gamma}

\newcommand{\thickness}{\vartheta}
\newcommand{\len}{s}

\newcommand{\Len}{S}

\newcommand{\pdlen}{\partial_{\len}}

\newcommand{\boper}{\theoperator{B}(\dspos,\dsdspos)}
\newcommand{\hoper}{\theoperator{H}(\dspos,\thickness)}

\newcommand{\kopermat}{\theoperator{K}_{\mathrm{material}}(\dspos,\dsdspos)}
\newcommand{\kopergeo}{\theoperator{K}_{\mathrm{geometric}}(\dspos,\dsdspos)}
\newcommand{\stress}{\thevector{s}}
\newcommand{\strain}{\thevector{e}}
\newcommand{\elasstress}{\stress_{\mathrm{elastic}}}
\newcommand{\elasstrain}{\strain_{\mathrm{elastic}}}

\newcommand{\viscstress}{\stress_{\mathrm{fluid}}}

\newcommand{\difflen}{d\len}

\newcommand{\diffinterface}{\mathrm{d}\gamma}

\newcommand{\coor}{\thevector{q}}

\newcommand{\varcoor}{\delta\coor}

\newcommand{\nmat}{\thematrix{N}}
\newcommand{\dsnmat}{\nmat'}
\newcommand{\dsdsnmat}{\nmat''}
\newcommand{\bmat}{\thematrix{B}}
\newcommand{\hmat}{\thematrix{H}}

\newcommand{\diffxi}{\mathrm{d}\xi}

\newcommand{\cgg}[1]{\textcolor{red}{CGG:#1}}

\begin{document}

\maketitle

\begin{abstract}
We propose a finite element method for simulating one-dimensional solid models moving and experiencing large deformations while immersed in generalized Newtonian fluids. The method is oriented towards applications involving microscopic devices or organisms in the soft-bio-matter realm. By considering that the strain energy of the solid may explicitly depend on time, we incorporate a mechanism for active response. The solids are modeled as Cosserat rods, a detailed formulation being provided for the special case of a planar non-shearable rod. The discretization adopts one-dimensional Hermite elements for the rod and low-order Lagrange two-dimensional elements for the fluid's velocity and pressure. The fluid mesh is boundary-fitted, with remeshing at each time step.  Several time marching schemes are studied, of which a semi-implicit scheme emerges as most effective. The method is demonstrated in very challenging examples: the roll-up of a rod to circular shape and later sudden release, the interaction of a soft rod with a fluid jet and the active self-locomotion of a sperm-like rod. The article includes a detailed description of a code that implements the method in the Firedrake library. 
\end{abstract}


% -*- mode: latex; TeX-master: "main.tex" -*-

\section{Introduction}

One-dimensional solid models (strings, cables, trusses, bars, beams, filaments, rods, etc.) are widely used in macroscopic structural analysis. Well established procedures exist, together with academic and commercial codes, that approximate these models and provide meaningful solutions to the analyst.
They are also useful in modeling microscopic phenomena, both in the solid realm (e.g., fiber-reinforced composites \cite{Le2017, Steinbrecher2020}) and in the soft-bio-matter realm (e.g., human hair \cite{Bertails2006, Kmoch2009}, climbing plants \cite{McMillen2002},  catheters and stents \cite{Wang2017}, DNA \cite{Benham1979, Schlick1995}, cilia and flagella \cite{RodrigoVelezCordero2012} and further biological applications \cite{Fauci2001}) that is the focus of this contribution. 

Simple elongated organisms, or the appendices (flagella, cilia) of more complex ones, are frequently described as one-dimensional entities. This is very convenient, since the model unknowns reduce, for each time instant, to a set of functions defined on some real interval $[0,S]$. Initial conditions and control actions are in this way much simpler to build and analyze. Though much of the theory and methods of macroscopic solid mechanics apply straightforwardly, some challenges are raised, or emphasized, by the microscopic realm. 

Firstly, the usual boundary conditions that attach or clamp the solid body to a rigid inmovable boundary are in general absent. The body is often floating in the surrounding fluid, and if the action of the fluid is neglected, {\em or if the forces that the fluid exerts on the solid are frozen in time} (as done in staggered FSI algorithms), the physical/mathematical problem for the solid becomes undetermined (only computable up to an arbitrary translation and rotation). This difficulty is typical of sedimenting and swimming bodies, even if they are zero-dimensional (i.e., a set of rigid particles, connected or not). Adequate treatments of the translation and rotation degrees of freedom (dofs) of bodies immersed in fluid have been presented in \cite{Maniyeri2012, Ishimoto2019, Paz2020} among others.

A second challenging aspect of modeling one-dimensional bodies in soft-bio matter is the very large displacements and deformations that are usually involved. The strong nonlinearities that appear in the differential geometry of curves comes into play and cannot be simplified/linearized. This makes that only quite sophisticated formulations of Computational Solid Mechanics (CSM) can be adopted. Later on we describe one such formulation, which corresponds to a planar nonlinear non-shearable rod, but it is certainly not the only possible one. The reader is referred to \cite{Romero2020, Gebhardt2021} for further details in the context of similar formulations in the three-dimensional setting.

Another peculiarity of soft-bio-matter problems is that, as happens also in aerodynamics, the effect of the surrounding fluid on the solid bodies is much more consequential than simply energy dissipation. Swimmers, in particular, rely on the complex response of the fluid for their propulsion and guidance \cite{Childress2012, Park2019, Danis2019, Liu2020}. A quite accurate fluid model is thus necessary to approximate the motion of microscopic organisms. 
In many cases, this means solving the full Navier-Stokes equations or variants thereof (in the case, for example, of non-Newtonian fluid's rheology \cite{Shen2011, MontenegroJohnson2013}), since simplified models such as Ressistive Force Theory may yield misleading conclusions \cite{Rorai2018}. 

In CSM it is customary to denote as Fluid-Structure Interaction (FSI) problems those that require the coupled solution of the fluid's and solid's partial differential equations. Thus, one of the challenges of soft-bio problems is that they are often strongly coupled FSI problems. The literature on FSI is quite vast and many previous works introduce concepts and methods that justify and motivate some of the procedures proposed in this article \cite{Farhat1998, Sarrate2001, Quaini2007, BoffiGastaldi2017, annese2020splitting}.

We would also like to mention that microscopic organisms are usually {\em active}, which implies moreover, the need of considering them as controllable or controlled systems \cite{Alouges2013a, Alouges2013b, Faris2015}. They produce mechanical energy through molecular motors that consume chemical compounds such as ATP ({\em adenosine triphosphate}). The incorporation of active elements in the model, which is not standard in CSM, is quite mandatory in the soft-bio realm.

In the following we propose finite element algorithms for one-dimensional solid models moving and deforming within a fluid. The algorithm does not rely on the regularizing effect of inertia of either the solid or the fluid, and will be presented neglecting all inertial terms to simplify the exposition and notation. The formulation adopts moving boundary-following meshes, which are in general more accurate than immersed-boundary methods \`a la Peskin \cite{Peskin2000} at the expense of requiring a robust meshing software. It should be mentioned that the complexity of micro-organisms' shape and movement is simpler to address using boundary element and other singularity-based methods. They have been developed and exploited in a quite extensive body of previous work (see \cite{PhanThien1987, Trachtenberg2003, Goto2005, Smith2009, Shum2015, Shum2019, Giuliani2021} and references therein), but they have limitations that make a finite element method {\em for the fluid} to be interesting. 
%In the context of microswimmers, those popular models are not robust enough, and when convergent may conduct to misleading conclusions \cite{Rorai2018}. 
The most evident limitation is that singularity-based methods cannot be applied to fluids with non-Newtonian rheology, which are quite abundant. The second limitation arises when the body motion is coupled to the transport of some species (e.g., nutrients) by the surrounding fluid. Singularity-based methods do not compute a sparse representation of the bulk velocity of the fluid, which thus requires a (costly) additional reconstruction step before the transport calculations.

\section{Formulating the problem}

\subsection{The surface motion of a Cosserat rod}

A quite general model of rod-like solids is provided by the Cosserat theory of rods \cite{Antman1995, Rubin2000}. In this theory, a rod is a one-dimensional object which can thus be parameterized by a single {\bf body coordinate} $s\,\in\,[0,S]$. A {\bf configuration} of the rod is defined by two elements: (a) a continuously differentiable mapping from $[0,S]$ into $\mathbb{R}^d$, where $d$ is the dimension of the ambient space, satisfying $m'(s)\neq 0$ for all $s$; and (b) a set of $d$ vector-valued functions defined in $[0,S]$, called {\bf directors}. The first director at each point is the tangent vector, and the other (or others) must be {\em linearly independent} so as to complete a basis of $\mathbb{R}^d$. 

The reference configuration $(\mathbf{M},\mathbbmss{D})$ of the rod is denoted by the mapping $\mathbf{M}(s)$ and the directors $\mathbbmss{D}(s)=\{\mathbf{D}_i(s)\}$, with $i=1,\ldots,d$ and $s\,\in\,[0,S]$. 
%Adopting a cartesian coordinate system in $\mathbb{R}^d$, fixed once and for all, we can identify the set of director fields $\mathbf{D}(s)$ with a $d\times d$ invertible matrix $D(s)$.
Let $\mathbf{m}$ and $\mathbbmss{d}=\{\mathbf{d}_i\}$, with $\mathbf{d}_1(s)=\mathbf{m}'(s)$, be the analogous functions that characterize another arbitrary configuration (the ``actual'' configuration). 
The basic idea is to adopt kinematical assumptions such that the deformation of the body's wet surface is totally encoded in the functions $\mathbf{m}$ and $\{\mathbf{d}_i\}$. Let $\widehat{\mathcal{S}}$ be the wet surface in the reference configuration, assumed to be closed and orientable. The material point that occupies the position $\mathbf{Y}\,\in\,\widehat{\mathcal{S}}$ in the reference configuration is assigned a body coordinate $\sigma(\mathbf{Y})\,\in\,[0,S]$. The function $\sigma:\widehat{\mathcal{S}}\to [0,S]$ is assumed continuous and differentiable. If $\widehat{\mathcal{C}}=\{\mathbf{M}(s), s\in [0,S]\}$ is the {\bf centroidal line} of the body in the reference configuration, $\mathbf{M}(\sigma(\mathbf{Y}))$ is the position of the centroidal point to which $\mathbf{Y}$ is considered ``attached''. We can thus write, for all $\mathbf{Y}\in\widehat{\mathcal{S}}$,
\begin{equation}
\mathbf{Y}=\mathbf{M}(\sigma(\mathbf{Y})) + \sum_{k=1}^d
\theta_k(\mathbf{Y})\,\mathbf{D}_k (\sigma(\mathbf{Y}))~,
\end{equation}
which uniquely defines $d$ real functions $\{\theta_i\}$ on $\widehat{\mathcal{S}}$. In matrix notation
\begin{equation}
 \boldsymbol{\theta}(\mathbf{Y})=\left [ \mathbbmss{D}(\sigma(\mathbf{Y}))\right ]^{-1}\,\left ( \mathbf{Y} - \mathbf{M}(\sigma(\mathbf{Y})\right )~, \label{eq:titasgen}
 \end{equation}
 where $\mathbbmss{D}$ is the matrix that has as columns the components of the directors $\{\mathbf{D}_i\}$ in some (Cartesian for simplicity) basis $\{\mathbf{E}_i\}$, $i=1,\ldots,d$.
Also, it is always possible (upon reparameterizing the body so that $s$ is the arc length along $\widehat{\mathcal{C}}$) to choose the reference directors forming an orthonormal basis, in which case $\mathbbmss{D}$ is an orthogonal matrix with determinant $+1$, i.e., an element of $SO(d)$.

The kinematical assumption is that in the actual configuration, the position $\mathbf{y}$ of the same material point is given by
\begin{equation}
\mathbf{y}=\Psi(\mathbf{Y})=\mathbf{m}(\sigma(\mathbf{Y}))+\sum_{k=1}^d \theta_k(\mathbf{Y})\,\mathbf{d}_k (\sigma(\mathbf{Y}))~.
\label{eq:ygeneral}
\end{equation}
This can again be rewritten in matrix notation as
\begin{equation}
\mathbf{y} = \mathbf{m}(\sigma(\mathbf{Y})) + \mathbbmss{d}\,
(\sigma(\mathbf{Y}))\,\boldsymbol{\theta}(\mathbf{Y})~,
\end{equation}
where $\mathbbmss{d}$ is the matrix that has as columns the components of the directors $\{\mathbf{d}_k\}$ in the basis $\{\mathbf{E}_i\}$, $i=1,\ldots,d$.

Notice that for each $s\in [0,S]$, the set $\widehat{\mathcal{R}}(s)=\sigma^{-1}(s)$ consists of the points in the reference wet surface $\widehat{\mathcal{S}}$ that have body coordinate $s$. From (\ref{eq:ygeneral}) we see that the actual position of this set of {\em material points} is totally dictated by $\mathbf{m}(s)$ and $\mathbbmss{d}(s)$. Further, let $\mathbf{A}$ and $\mathbf{B}$ be two points belonging to $\widehat{\mathcal{R}}(s)$. Then
\begin{equation}
\mathbf{b}-\mathbf{a}=\mathbbmss{d}(s)\,\left [ \boldsymbol{\theta}(\mathbf{B}) - \boldsymbol{\theta} (\mathbf{A}) \right ]
= \mathbbmss{d}(s)\,[\mathbbmss{D}(s)]^{-1}
\left ( \mathbf{B} - \mathbf{A}\right )~,
\end{equation}
so that, restricted to $\widehat{\mathcal{R}}(s)$, the deformation is {\em linear}. Later on we consider bodies with a head attached to $s=S$, in the sense that $\widehat{\mathcal{R}}(S)$ is a region of positive measure of $\widehat{\mathcal{S}}$ (see Fig \ref{fig:refswimmer}). Such a head is thus subject to the translation $\mathbf{m}(S)-\mathbf{M}(S)$ plus the linear deformation of gradient $\mathbbmss{d}(S)[\mathbbmss{D}(S)]^{-1}$. A similar appendage can be attached to $s=0$.

Assuming that $\mathbf{m}$, $\{\mathbf{d}_k\}$, $\sigma$ and $\{\theta_k\}$ are regular enough, the mapping $\Psi$ is continuous and differentiable. Under the additional restriction of $\Psi$ being injective (no self intersections) we have that the actual wet surface $\mathcal{S}=\Psi(\widehat{\mathcal{S}})$ is also closed and orientable because $\Psi:\widehat{\mathcal{S}}\to \mathcal{S}$ is a homeomorphism.

In this work we do not attempt a rigorous characterization of the configurational manifold $\mathcal{Q}$ of a Cosserat rod. It is certainly a subset of the product space $V \times W$, where $V$ is the space of differentiable vector-valued functions from $[0,S]$ to $\mathbb{R}^d$ and $W$ is the space of matrix-valued functions from $[0,S]$ to $\mathbb{R}^{d\times d}$. The differentiability of $\mathbf{m}$ and $\mathbbmss{d}$ must ensure not only that $\Psi$ is a homeomorphism, but also that the mechanical energy with which the rod is endowed is finite. The latter is frequently the most restrictive condition.

The tangent space $T_{\mathbf{q}}\mathcal{Q}$ of configurational variations $\delta \mathbf{q}=(\delta \mathbf{m},\delta \mathbbmss{d})$ with respect to the actual configuration $\mathbf{q}=(\mathbf{m},\mathbbmss{d})$, can thus be identified to a {\em subspace} of $V\times W$, subject to the admissibility constraint $\delta \mathbf{d}_1=\delta \mathbf{m}'$.

Each configurational variation $\delta \mathbf{q}$ induces, through (\ref{eq:ygeneral}), a vector field $\mathbf{w}_{\delta \mathbf{q}}$ over the actual wet surface $\mathcal{S}$ given by
\begin{equation}
\mathbf{w}_{\delta \mathbf{q}}(\mathbf{y})=\delta \mathbf{y} = \delta \mathbf{m} + \delta \mathbbmss{d}\,\boldsymbol{\theta}~.
\label{eq:wdeltaq}
\end{equation}
It represents the infinitesimal movement of the material points of the surface produced by $\delta \mathbf{q}$. The operator defined by $H_\mathbf{q}:\delta \mathbf{q} \mapsto \mathbf{w}_{\delta \mathbf{q}}(\cdot)$
{\em is linear} from $T_{\mathbf{q}}\mathcal{Q}$ into $U(\mathcal{S})$, the space of vector fields on $\mathcal{S}$. 
%We will denote the image of $H_\mathbf{q}$ as $K(\mathcal{S})$, i.e.,
%\begin{equation}
%K(\mathcal{S}(\mathbf{q}))=\left \{
%\mathbf{w}: \mathcal{S}(\mathbf{q}) \to \mathbb{R}^d, \mbox{ such that } \exists \mathbf{p}\in T_{\mathbf{q}}\mathcal{Q} \mbox{ satisfying }
%\mathbf{w}=H_{\mathbf{q}}\mathbf{p} \right \}.
%\end{equation}
%$K(\mathcal{S}(\mathbf{q}))$ consists of all possible surface velocity fields of the rod model when its configuration is $\mathbf{q}$.

Letting now the actual configuration to evolve in time, dependence on the additional time variable $t$ appears in the expressions above. For functions of both $s$ and $t$ we keep the notation $f'$ for $\partial_s f$ and introduce $\dot{f}$ for $\partial_t f$.
  
It is clear that the configuration $\mathbf{q}(t)$ and the generalized velocity $\mathbf{\dot{q}}(t)$ of the rod at a given instant $t$ must belong to $\mathcal{Q} \times T_{\mathbf{q}(t)}\mathcal{Q}$. From the definitions above, if $\mathbf{u}_b(\mathbf{x},t)$ is the velocity field of the body at the wet surface, then it is computed from $\mathbf{\dot{q}}$ as
\begin{equation}
\mathbf{u}_b(\mathbf{y},t) = \mathbf{w}_{\mathbf{\dot{q}}(t)}(\mathbf{y},t)=[H_{\mathbf{q}(t)}\,\mathbf{\dot{q}}(t)](\mathbf{y},t)~,\qquad \forall \mathbf{y}\,\in\, \mathcal{S}(t)~.
\label{eq:ubq}
\end{equation}

\subsection{The case of a planar non-shearable rod} \label{sectrodkin}

Further kinematical assumptions can be introduced, leading to more specific models. They will result in further linear admissibility constraints in $T_{\mathbf{q}} \mathcal{Q}$. Let us consider a planar ($d=2$) non-shearable rod, which is adopted hereafter as application example. This rod model assumes that, in all configurations, $\mathbf{d}_2$ is always normal to $\mathbf{d}_1$, so that
\begin{equation}
\mathbf{d}_2 (s) = a(s)\,\mathbf{\check{n}}(s)=a(s)\,\mathbf{J}\,\mathbf{\check{t}}(s)=\frac{a(s)}{\|\mathbf{m}'(s)\|}\,\mathbf{J}\,\mathbf{m}'(s)~,
\end{equation}
where $\mathbf{\check{t}}$ and $\mathbf{\check{n}}$ are tangent and normal unit vectors, respectively, and
\begin{equation}
	\jye := 
	\begin{pmatrix}
		 0 &-1\\
    	 1 & 0
    \end{pmatrix} . \label{eq:Jota}
\end{equation}
We will further assume that the coefficient $a(s)$ is totally determined by the curve $\mathbf{m}$. An example of this situation is when $a(s)$ is a function of $\mathbf{m}'(s)$, $\mathbf{m}''(s)$, etc., but the dependence could as well be nonlocal. To fix ideas, let us consider the particular case 
\begin{equation}
 a=a\left (\|\mathbf{m}'(s)\|\right )~. \label{eq:aofsgen}
\end{equation} 
With these kinematical assumptions the function $\mathbf{m}(\cdot)$ totally determines both directors and thus also the whole deformation of the wet surface, as given by
 \begin{equation} \mathbf{y}=\mathbf{m}(s)+\theta_1(\mathbf{Y})\,\mathbf{m}'(s)+\theta_2(\mathbf{Y})\,\frac{a(\|\mathbf{m}'(s)\|)}{\|\mathbf{m}'(s)\|}\,\jye\,\mathbf{m}'(s)
 \end{equation}
 with $s=\sigma(\mathbf{Y})$ and $\mathbf{y}=\Psi(\mathbf{Y})$. The additional constraint on $T_{\mathbf{q}}\mathcal{Q}$ is, defining $b(s)=\|\mathbf{m}'(s)\|$,
 \begin{equation}
 \delta \mathbf{d}_2 (s) = \frac{a'(b(s))}{b(s)} \mathbf{m}'(s)\cdot \delta \mathbf{m}'(s)\,\mathbf{\check{n}}(s)+ a(b(s))\,\delta \mathbf{\check{n}}(s)~.
 \end{equation}
 Since
 \begin{equation}
 \delta \mathbf{\check{n}}(s) = \jye\,\delta \mathbf{\check{t}}(s),\qquad
 \delta \mathbf{\check{t}}(s) = \frac{1}{b(s)} \mathbf{P}^{\perp}\,\delta \mathbf{m}'~,
 \end{equation}
 where $\mathbf{P}^{\perp}=\mathbf{I}-\mathbf{\check{t}}\otimes\mathbf{\check{t}}$, we can also write the constraint as
 \begin{equation}
 \delta \mathbf{d}_2 (s) = \left [\frac{a'(b(s))}{b(s)}  \mathbf{\check{n}}(s)\otimes \mathbf{m}'(s)\,+ \frac{a(b(s))}{b(s)}\,\jye\,\mathbf{P}^{\perp}(s) \right ]\,
 \delta \mathbf{m}'(s)~.
 \label{eq:deltad2}
 \end{equation}
 
 The previous calculations show that (\ref{eq:wdeltaq}) in this case takes que form
 \begin{eqnarray}
 \mathbf{w}_{\delta \mathbf{q}}(\mathbf{y}) &=& \delta \mathbf{m} (s) + \left [ \theta_1(\mathbf{Y})\, \mathbf{I} + \theta_2(\mathbf{Y})\,\mathbf{G}(s)\right ]\,\delta \mathbf{m}'(s) \nonumber \\ & = &
 \left \{ \mathbf{I} + \left [ \theta_1(\mathbf{Y})\, \mathbf{I} + \theta_2(\mathbf{Y})\,\mathbf{G}(s)\right ]\,\partial_s \right \}~\delta \mathbf{m}(s)~,
 \end{eqnarray}
 where $\mathbf{G}(s)$ is the $2\times 2$ matrix that appears between brackets on the right-hand side of (\ref{eq:deltad2}).

Without loss of generality, we assume that $\|\mathbf{D}_1(s)\|=\|\mathbf{D}_2(s)\|=1$, so that the reference directors coincide with tangent and normal unit vectors,
\begin{equation}
\mathbf{D}_1(s)=\mathbf{\check{T}}(s),\qquad
\mathbf{D}_2(s)=\mathbf{\check{N}}(s)~, \qquad \forall\,s.
\end{equation}
In such a case, the coefficient $a(s)$ measures the local contraction (if $<1$) or expansion (if $>1$) of imaginary fibers along the normal direction, while $b(s)=\|\mathbf{m}'(s)\|$ is the local contraction (if $<1$) or expansion (if $>1$) of imaginary fibers along the tangent direction.

The configuration of this model of rod, as becomes clear from above, is totally determined by $\mathbf{m}(\cdot)$. One can thus adopt a {\em reduced} configuration manifold $\mathcal{Q}$ consisting just of the centroid mappings $\mathbf{m}:[0,S]\to \mathbb{R}^d$, in which case $\mathbf{q}=\mathbf{m}$, $\delta \mathbf{q}=\delta \mathbf{m}$ (unrestricted), and $T_\mathbf{q}\mathcal{Q}$ is a subspace of $V$. In this reduced formulation, the operator $H_{\mathbf{m}}:T_{\mathbf{q}}\mathcal{Q}\to U(\mathcal{S})$ (the space of vector fields on $\mathcal{S}$) is given by
\begin{equation}
H_\mathbf{m} = \mathbf{I} + \left [ \theta_1(\mathbf{Y})\, \mathbf{I} + \theta_2(\mathbf{Y})\,\mathbf{G}(s)\right ]\,\partial_s ~.
\end{equation} 
Furthermore, expression (\ref{eq:ubq}) for the velocity of the body at the wetted surface particularizes to
\begin{equation}
\mathbf{u}_b(\mathbf{y},t)=[H_{\mathbf{m}(\cdot,t)}\mathbf{\dot{m}}(\cdot,t)]\,(\mathbf{y},t)
= \mathbf{\dot{m}}(s,t)+ \left [ \theta_1(\mathbf{Y})\, \mathbf{I} + \theta_2(\mathbf{Y})\,\mathbf{G}(s) \right ] \, \mathbf{\dot{m}}'(s,t)~. \label{eq:inducedvel}
\end{equation}

\subsection{The ambient fluid and the fluid-solid kinematics}

Let us denote by $\Omega_f(\mathbf{q})$ the domain of the ambient fluid when the rod is in the configuration $\mathbf{q}$. It has $\mathcal{S}(\mathbf{q})$ as internal boundary and let $\Gamma$ be its (fixed) external boundary of which a portion $\Gamma_{\mathbf{u}}$ is a wall at which the fluid velocity $\mathbf{u}$ is known (and equal to zero, for simplicity). Considering just the fluid part of the system, the space of admissible velocity fields is then
\begin{equation}
X_{\mathbf{q}}=\left \{ \mathbf{w}\,\in\,\left [ H^1(\Omega_f(\mathbf{q})) \right ]^d~|~\mathbf{w}(\mathbf{x})=\mathbf{0}~\mbox{ a.e. in } \Gamma_{\mathbf{u}} \right \}~.
\end{equation}

{\bf Remark:} The space $H^1(\Omega_f)$ in the definition of $X_{\mathbf{q}}$ corresponds to a Newtonian ambient fluid. Other rheological behaviors may require its substitution by a suitable Sobolev space $W^{m,p}(\Omega_f)$. 

The most frequent kinematical conditions at the fluid-solid interface are {\em no-slip conditions}, i.e.,
 the fluid velocity $\mathbf{u}$ must satisfy
 \begin{equation}
  \mathbf{u}(\mathbf{y})=\mathbf{u}_b(\mathbf{y})=\left [ H_{\mathbf{q}}\,\mathbf{\dot{q}}\right ] (\mathbf{y})~,\qquad \qquad \forall\,\mathbf{y}\,\in\,\mathcal{S}~.
  \end{equation}
  
 We assume that the kinematics of the solid at configuration $\mathbf{q}$ is sufficiently smooth so that for any $\mathbf{p}\,\in\,T_{\mathbf{q}}\mathcal{Q}$ the induced surface field $\mathbf{w}=H_{\mathbf{q}}\mathbf{p}$ can be extended to $\Omega_f$ as an element of $X_{\mathbf{q}}$, satisfying $\|\mathbf{w}\|_{X_{\mathbf{q}}} \leq c\, \|\mathbf{p}\|_{T_{\mathbf{q}}\mathcal{Q}}$ for some $c\geq 0$. This requires, in particular, that $H_{\mathbf{q}}\mathbf{p}$ belongs to the trace space of $X_{\mathbf{q}}$ on $\mathcal{S}(\mathbf{q}))$ for all $\mathbf{p}$. At this point it is convenient to define
   \begin{equation}
   \widehat{H}_{\mathbf{q}}=L_{\mathbf{q}}H_{\mathbf{q}}
   \label{eqdefHhat}
 \end{equation}
 as a linear continuous operator from $T_{\mathbf{q}}\mathcal{Q}$ into $X_{\mathbf{q}}$, where $L_{\mathbf{q}}$ is an (arbitrary) continuous linear extension operator from functions defined on $\mathcal{S}(\mathbf{q})$ to functions in $X_{\mathbf{q}}$. 
  
  As a consequence, at any instant $t$ the pair $(\mathbf{\dot{q}},\mathbf{u})\,\in\,T_{\mathbf{q}(t)}\mathcal{Q}\times X_{\mathbf{q}(t)}$, which we call the {\em generalized velocity of the fluid-solid system}, must belong to the {\em space of kinematically admissible fluid-solid motions} defined by
 \begin{equation}
Z_{\mathbf{q}(t)}=\left \{
(\mathbf{p},\mathbf{w})\,\in\,T_{\mathbf{q}(t)}\mathcal{Q}\times
X_{\mathbf{q}(t)} ~|~\left. \mathbf{w}\right |_{\mathcal{S}(\mathbf{q}(t))} = H_{\mathbf{q}(t)}\,\mathbf{p}\right \}~.
\end{equation}
where the notation $\left.f\right |_{\mathcal{S}}$ denotes the restriction of the function to $\mathcal{S}$ (in the sense of traces). The variational formulation of the instantaneous FSI problem at time $t$ is posed in the space $Z_{\mathbf{q}(t)}$.

\subsection{The principle of virtual work and the evolution problem}

In the previous sections we have defined and described the kinematics of the fluid-solid systems of which the solid's motion is parameterized as a Cosserat rod. Now we consider the dynamics of such systems, as dictated by the principle of virtual work.

We assume that the rod is endowed with a differentiable energy $E(t,\mathbf{q})$, dependent on both the time and the configuration, and that all the inertial effects are neglected. Then, at each instant $t$, the instantaneous configuration $\mathbf{q}\,\in\,\mathcal{Q}$ and generalized velocity $(\mathbf{\dot{q}},\mathbf{u})\,\in\,Z_{\mathbf{q}}$ must satisfy
\begin{equation}
\langle D_{\mathbf{q}}E(t,\mathbf{q}) ,  \delta \mathbf{q}\rangle  + \int_{\Omega_f} \boldsymbol{\sigma}\bm{:} \boldsymbol{\varepsilon}(\delta {\mathbf w})~d\Omega
= \int_{\Omega_f} \mathbf{f}\bm{\cdot} \delta \mathbf{w}~d\Omega + \langle \mathbf{F}(t), \delta \mathbf{q}\rangle
\label{eqpvw}
\end{equation}
for all $(\delta \mathbf{q},\delta \mathbf{w})\,\in\,Z_{\mathbf{q}}$, where $D_{\mathbf{q}}E(\mathbf{q})$ is the energy differential (with respect to $\mathbf{q}$), $\boldsymbol{\sigma}$ is the Cauchy stress tensor of the fluid, $\boldsymbol{\varepsilon}$ is the symmetric gradient operator, $\mathbf{f}$ is the body force acting on the fluid, $\mathbf{F}$ is the generalized (non-conservative) force acting on the rod not coming from the ambient fluid and $\langle \cdot,\cdot \rangle$ is the duality pairing between $T_{\mathbf{q}}\mathcal{Q}$ and its dual. In what follows we will assume that $\mathbf{f}=0$ (notice that the Dirichlet boundary conditions are also zero). It is not at all essential, but it simplifies the presentation considerably.

From (\ref{eqpvw}) one obtains the evolution problem in variational form. Let us consider an incompressible Stokesian fluid model defined by
\begin{equation}
\boldsymbol{\sigma}=-p\,\mathbf{I} + 2\,\mu\,\boldsymbol{\varepsilon}(\mathbf{u})~, \qquad
\nabla \bm{\cdot} \mathbf{u}=0
\end{equation}
(incompressible Newtonian fluid, or quasi-Newtonian if $\mu$ depends on the deformation rate), where $p$ stands for the pressure field. Then the evolution problem reads

\bigskip

\noindent {\bf Problem P1:} ``Find $\mathbf{q}:(0,T)\to \mathcal{Q}$, such that for a.e. $t\,\in\,(0,T)$ there exists $\mathbf{u}(\cdot,t)\,\in\,X_{\mathbf{q}(t)}$ and $p(\cdot,t)\,\in\,\Pi_{\mathbf{q}(t)}=L^2(\Omega_f(\mathbf{q}(t)))$ satisfying $\mathbf{q}(0)=\mathbf{q}_0\,\in\,\mathcal{Q}$, together with
\begin{equation}
  (\mathbf{\dot{q}(t)},\mathbf{u}(\bm{\cdot},t))\,\in\,Z_{\mathbf{q}(t)}\qquad
  \mbox{  (i.e., kinematic admissibility),}
  \label{eq18}
\end{equation}
and 
 \begin{eqnarray}
\langle D_{\mathbf{q}}E(t,\mathbf{q}(t)) ,  \delta \mathbf{q}\rangle  + \int_{\Omega_f} 2\mu\, \boldsymbol{\varepsilon}(\mathbf{u})\bm{:} \boldsymbol{\varepsilon}(\delta {\mathbf w})~d\Omega - \int_{\Omega_f} p \,\nabla \bm{\cdot}\,\delta \mathbf{w} ~d\Omega
&=& %\int_{\Omega_f} \mathbf{f}\bm{\cdot} \delta \mathbf{w}~d\Omega +
\langle \mathbf{F}, \delta \mathbf{q}\rangle ,\label{eqp1a}\\
\int_{\Omega_f} r\,\nabla \bm{\cdot} \mathbf{u} ~d\Omega &=& 0~, \label{eqp1b}
\end{eqnarray}
 for all $(\delta \mathbf{q},\delta \mathbf{w})\,\in\,Z_{\mathbf{q}(t)}$ and all $r\,\in\,\Pi_{\mathbf{q}(t)}$''.

 Under realistic hypotheses on the data, problem P1 is expected to admit a unique solution. These hypotheses involve the energy function, requiring it to be smooth enough and exhibiting some kind of suitable convexity, with the kernel of $D_{\mathbf{q}}E$ (typically rigid-body motions) of finite dimension. They also involve the forcing, requiring it to be bounded in some suitable norm, for example 
 \begin{equation}
   \sup_{t\,\in\,(0,T)} \|\mathbf{F}(t)\|_{T_{\mathbf{q}(t)}Q'} ~< ~+\infty~;
 \end{equation}
 and the geometry of the system, requiring the operator $\widehat{H}_{\mathbf{q}}$ to be {\bf injective} and {\bf compact} from $T_{\mathbf{q}}\mathcal{Q}$ into $X_{\mathbf{q}}$. In particular, any non-zero infinitesimal deformation $\delta \mathbf{q}$ of the solid must produce a non-zero surface field $H_{\mathbf{q}}\delta \mathbf{q}$ on $\mathcal{S}$. In particular, if the rod is not attached to any solid boundary, the rigid motions are controlled by the viscous operator.

 %In fact, a more technical hypothesis is needed here to state the result with some rigor. Remember the operator $\widehat{H}_{\mathbf{q}}$ which to each $\mathbf{p}\in T_{\mathbf{q}}\mathcal{Q}$ assigns a velocity field $\widehat{H}_{\mathbf{q}}{\mathbf{p}}\,\in\,X_{\mathbf{q}}$ over $\Omega_f$. It is necessary that given any bounded sequence $\mathbf{p}_1,\mathbf{p}_2,\ldots $ in $T_{\mathbf{q}}\mathcal{Q}$ the sequence of images  $\widehat{H}_{\mathbf{q}}\mathbf{p}_1,\widehat{H}_{\mathbf{q}}\mathbf{p}_2,\ldots $ has at least one convergent subsequence in $X_{\mathbf{q}}$. In other words, the operator  $\widehat{H}_{\mathbf{q}}$ needs to be compact.

 The numerical approximation proposed in this work is best understood if we introduce a new problem, equivalent to P1. For this purpose, let us consider an arbitrary instant $t$ at which the solid configuration is $\mathbf{q}$. Then, for each $\mathbf{p}\,\in\,T_{\mathbf{q}}\mathcal{Q}$ there exists a unique velocity-pressure pair $(\mathbf{u},p)$ that satisfies the boundary condition
 \begin{equation}
   \mathbf{u}|_{\mathcal{S}(\mathbf{q})}=H_{\mathbf{q}}\mathbf{p}
 \end{equation}
 and the fluid's differential equations in $\Omega_f(\mathbf{q})$. In variational terms, it is the unique solution to (\ref{eqp1a})-(\ref{eqp1b}) when $\delta \mathbf{q}=0$, that is,
 \begin{eqnarray}
 \int_{\Omega_f(\mathbf{q})} 2\mu\, \boldsymbol{\varepsilon}(\mathbf{u})\bm{:} \boldsymbol{\varepsilon}(\delta {\mathbf w})~d\Omega - \int_{\Omega_f(\mathbf{q})} p \,\nabla \bm{\cdot}\,\delta \mathbf{w} ~d\Omega
 &=& 0 %\int_{\Omega_f} \mathbf{f}\bm{\cdot} \delta \mathbf{w}~d\Omega
 ~,\label{eqp1a2}\\
\int_{\Omega_f(\mathbf{q})} r\,\nabla \bm{\cdot} \mathbf{u} ~d\Omega &=& 0~, \label{eqp1b2}
 \end{eqnarray}
 for all $r\,\in\,\Pi_{\mathbf{q}}$ and all $\delta \mathbf{w}$ in
\begin{equation}
X^0_{\mathbf{q}} = \left \{ \mathbf{v}\,\in\,X_{\mathbf{q}} ~|~ \mathbf{v}=0 \mbox{ on }\mathcal{S}(\mathbf{q}) \right \}~.
\end{equation}
We denote this velocity field by $\mathcal{U}_{\mathbf{q},\mathbf{p}}$ and the corresponding pressure field by $\mathcal{P}_{\mathbf{q},\mathbf{p}}$. By construction, given $\mathbf{q}$ and $\mathbf{p}$ it holds that
\begin{equation}
(\mathbf{p},\mathcal{U}_{\mathbf{q},\mathbf{p}})\,\in\,Z_{\mathbf{q}}
\end{equation}
and that
\begin{equation}
- \nabla \cdot \left [ 2 \mu\,\boldsymbol{\varepsilon}(\mathcal{U}_{\mathbf{q},\mathbf{p}}) \right ]
+ \nabla \mathcal{P}_{\mathbf{q},\mathbf{p}} = 0,\qquad
\nabla \cdot \mathcal{U}_{\mathbf{q},\mathbf{p}} = 0,~\qquad\qquad \mbox{ a.e. in }
\Omega_f(\mathbf{q})~.
\end{equation}
From the assumed continuity and injectivity of $\widehat{H}_{\mathbf{q}}$, we have that $\mathcal{U}_{\mathbf{q},\mathbf{p}}$ and $\mathcal{P}_{\mathbf{q},\mathbf{p}}$ are continuous in $\mathbf{p}$. If the fluid is Newtonian, they are {\em linear} in $\mathbf{p}$ and there exists $C_S\geq 0$ (dependent on $\mathbf{q}$) such that, for all $\mathbf{p}\neq 0$,
\begin{equation}
 0<\|\left ( \mathcal{U}_{\mathbf{q},\mathbf{p}},\mathcal{P}_{\mathbf{q},\mathbf{p}}\right )\|_{X_q\times \Pi_q} \leq C_S\,\|\mathbf{p}\|_{T_{\mathbf{q}}\mathcal{Q}}~.
\label{equbounded}
\end{equation}

The formulation of problem P1 that is considered for approximation is then as follows:

\bigskip

\noindent{\bf Problem P2:} ``Find $\mathbf{q}:(0,T)\to \mathcal{Q}$ such that $\mathbf{q}(0)=\mathbf{q}_0$ and
%\begin{equation}
%  \int_0^T \left [ b(\mathbf{q}(t);\mathbf{\dot{q}}(t),{\delta \mathbf{q}})
%    + \langle D_{\mathbf{q}}E(t,\mathbf{q}(t)),\delta \mathbf{q} \rangle \right ] ~dt
%    = \int_0^T \langle \mathbf{F}(t), \delta \mathbf{q} \rangle~dt
%\end{equation}
\begin{equation}
  b(\mathbf{q}(t);\mathbf{\dot{q}}(t),{\delta \mathbf{q}})
    + \langle D_{\mathbf{q}}E(t,\mathbf{q}(t)),\delta \mathbf{q} \rangle 
    =  \langle \mathbf{F}(t), \delta \mathbf{q} \rangle \label{eq:P2cont}
\end{equation}
for a.e. $t\,\in\,(0,T)$ and all $\delta \mathbf{q}$ such that $\delta \mathbf{q} (t)\,\in\,T_{\mathbf{q}(t)}\mathcal{Q}$,
where $b(\mathbf{q};\cdot,\cdot):T_{\mathbf{q}}\mathcal{Q}\times T_{\mathbf{q}}\mathcal{Q} \to
    \mathbb{R}$ is given by
\begin{equation}
  b(\mathbf{q};\mathbf{p},\delta \mathbf{q})=
  \int_{\Omega_f(\mathbf{q})} \left [ 2\mu\,\boldsymbol{\varepsilon} (\mathcal{U}_{\mathbf{q},\mathbf{p}}) \bm{:} \boldsymbol{\varepsilon} (\widehat{H}_{\mathbf{q}}\delta \mathbf{q})
-\mathcal{P}_{\mathbf{q},\mathbf{p}}\,\nabla \cdot (\widehat{H}_{\mathbf{q}}\delta \mathbf{q})
   \right ] ~d\Omega~.\mbox{''}
\end{equation}
Notice that $b(\mathbf{q};\mathbf{p},\delta \mathbf{q})$ may be nonlinear in $\mathbf{p}$ and is always linear in $\delta \mathbf{q}$.

\bigskip

\noindent{\bf Remark:} The fluid dissipation is given by
\begin{equation}
  \mathcal{D}(\mathbf{q};\mathbf{p})=\int_{\Omega_f(\mathbf{q})}2\mu \,\boldsymbol{\varepsilon} (\mathcal{U}_{\mathbf{q},\mathbf{p}}) \bm{:} \boldsymbol{\varepsilon} (\mathcal{U}_{\mathbf{q},\mathbf{p}})~d\Omega~,
\end{equation}
which using (\ref{eqp1a2})-(\ref{eqp1b2}) is seen to satisfy, for $\mathbf{p}\neq 0$, 
\begin{equation}
  \mathcal{D}(\mathbf{q};\mathbf{p})=b(\mathbf{q};\mathbf{p},\mathbf{p})>0~.  
\end{equation}

\bigskip

\noindent{\bf Remark:} In the linear (Newtonian) case $b(\mathbf{q};\cdot,\cdot)$ is a scalar product in $T_{\mathbf{q}}\mathcal{Q}$. In fact, the extension operator $L_{\mathbf{q}}$ (see (\ref{eqdefHhat})) can then be chosen as the Stokes solution; i.e. $L_{\mathbf{q}}H_{\mathbf{q}}\mathbf{p}=\mathcal{U}_{\mathbf{q},\mathbf{p}}$. Then the form $b(\mathbf{q};\cdot,\cdot)$ simplifies to
\begin{equation}
  b(\mathbf{q};\mathbf{p},\delta \mathbf{q})=
  \int_{\Omega_f(\mathbf{q})} 2\mu\,\boldsymbol{\varepsilon} (\mathcal{U}_{\mathbf{q},\mathbf{p}}) \bm{:} \boldsymbol{\varepsilon} (\mathcal{U}_{\mathbf{q},\delta\mathbf{q}})~d\Omega,
\end{equation}
which is indeed bilinear, symmetric and positive definite. This however will not be further exploited because our interest is in formulations that remain valid when the fluid is nonlinear.

\subsection{Full model of a planar non-shearable rod}

Let us now fully define problem P2 for the planar non-shearable rod of section \ref{sectrodkin}, so that the whole procedure becomes clear. We adopt the reduced configuration in which $\mathbf{q}=\mathbf{m}$, so that a configuration is a function $\mathbf{m}:[0,S]\to \mathbb{R}^2$. Similarly, the reference configuration is $\mathbf{M}:[0,S]\to \mathbb{R}^2$. 
It only remains to specify the energy function $E$. Assuming both $\mathbf{m}$ and $\mathbf{M}$ to be smooth enough, we define the rod's energy as
\begin{equation}
    E=\frac12 \int_0^S \left [ C_{\epsilon} \epsilon^2 + C_{\kappa} \kappa^2 \right ] ~ds~,
\end{equation}
where $\epsilon$ is the local measure of elongation change
\begin{equation}
  \epsilon = \mathbf{\check{t}}\cdot \mathbf{m}'-\mathbf{\check{T}}\cdot \mathbf{M}' - \epsilon_0~, \label{eq:epsilon}
\end{equation}
$\kappa$ measures the curvature change
\begin{equation}
  \kappa = \mathbf{\check{t}}\cdot \mathbf{\check{n}}'-\mathbf{\check{T}}\cdot \mathbf{\check{N}}'-\kappa_0~ \label{eq:kappa}
\end{equation}
and $C_{\epsilon}$ and $C_{\kappa}$ are positive constants (they could also depend on $s$). 
The functions $\epsilon_0(t,s)$ and $\kappa_0(t,s)$, assumed known, correspond to spontaneous values of elongation and curvature. The energy depends explicitly on time through them, thus providing a mechanism for the injection of mechanical energy into the system even if the external forces $\mathbf{F}$ are zero (for self-locomotion, for example). Broadly speaking, for $E(t,\mathbf{q}(t))$ to be finite, and thus $\mathbf{q}(t)\,\in\,\mathcal{Q}$, the tangential deformation must belong to $H^1(0,S)$ and the normal deformation to $H^2(0,S)$. This is sufficient regularity to expect that the induced vector fields $H_{\mathbf{q}(t)}\delta \mathbf{q}$ belong to $H^{\frac12}(\mathcal{S})$ and that $\widehat{H}_{\mathbf{q}}$ is injective and compact as required. In physical terms, this means that the smoothness of $\mathbf{m}(t,\cdot):[0,S]\to \mathbb{R}^2$ is dictated by the rod and not by the ambient fluid.

The first variation of the elastic energy reads
\begin{equation}
\delta E =  \langle D_{\mathbf{m}}E , \delta \mathbf{m} \rangle =
\int_0^S \left (\epsilon\,C_{\epsilon}\,\delta{\epsilon} + \kappa\,C_{\kappa}\,\delta{\kappa} \right ) ~\difflen 
\end{equation}
where the variations of $\epsilon$ and $\kappa$ are
\begin{eqnarray}
\delta{\epsilon} & = & \tanv \cdot \delta{\mathbf{m}} \label{eq:depsilon}\\
\delta{\kappa} & = & \moddspos^{-2}\dsdspos\pairing(2\norv\otimes\tanv-\jye)\,\delta{\mathbf{m}} - \moddspos^{-1}\norv \cdot \delta{\mathbf{m}} \label{eq:dkappa}
\end{eqnarray}
that can be written in compact form as
\begin{equation}
\delta E = 
  \int_0^S \left ( \boper\,\delta \mathbf{m} \right ) \cdot \elasstress~\difflen 
\end{equation}
with the stress resultants defined as
\begin{equation}
    \elasstress = \bm{C}\,\elasstrain,~
\end{equation}
where
\begin{equation}
\bm{C} = 
\begin{pmatrix}
C_{\epsilon} & 0 \\
0 & C_{\kappa}
\end{pmatrix}~,\quad
\elasstrain=
    \begin{bmatrix}
    \epsilon\\
    \kappa
    \end{bmatrix} . 
 %\checkmark
\end{equation}
The first component of $\elasstress$ is to be assigned to the elastic axial force and the second one to the elastic bending moment. The linearized strain operator is given by
\begin{equation}
    \boper := 
    \begin{pmatrix}
    \tanv\pairing\eye & 0 \\
    \moddspos^{-2}\dsdspos\pairing(2\norv\otimes\tanv-\jye) & \moddspos^{-1}\tanv\pairing\jye
    \end{pmatrix}
    \begin{pmatrix}
    \pdlen^1\\
    \pdlen^2
    \end{pmatrix} \label{eq:linearized}
\end{equation}
with $\pdlen^1:\pos\mapsto\dspos$ and $\pdlen^2:\pos\mapsto\dsdspos$ denoting differentiation operators.

Lastly, the second variation of the potential energy reads
\begin{equation}
     \Delta\delta E = \int\limits_0^\Len \varpos\pairing\left [ \kopermat+\kopergeo \right ] \Delta\pos~\difflen ,
\end{equation}
in which the operator associated to the material stiffness is given by
\begin{equation}
    \kopermat = \boper^T\bm{C}\boper
\end{equation}
and the operator associated to the geometrical stiffness is given by
\begin{equation}
    \kopergeo = \pd{\dspos}(\boper^T\elasstress)\pdlen^1+\pd{\dsdspos}(\boper^T\elasstress)\pdlen^2 \label{eq:Kgeom}.
\end{equation}

\begin{remark}
The formulation above is \textit{objective}, \ie invariant under translations and rotations, as well as \textit{path-independent}, \ie for conservative loading and time-independent energy function the work produced through any arbitrary closed path is identically zero. The rigid body motions span the null space of the stiffness operator.
\end{remark}

%%% Local Variables:
%%% mode: latex
%%% TeX-master: "main"
%%% End:

\section{Numerical method}

\subsection{Fully discrete formulation}

We make a strong hypothesis here, namely that finite element {\em spaces} $\mathcal{Q}_H\subset \mathcal{Q}$ exist such that any $\mathbf{q}\,\in\,\mathcal{Q}$ can be approximated by some $\mathbf{q}_H\,\in\,\mathcal{Q}_H$. Since $\mathcal{Q}_H$ is a vector space, it is identical to its tangent space at any point so that we will define, for simplicity,
\begin{equation}
  Q_H := \mathcal{Q}_H \equiv T_{\mathbf{q}_H}\mathcal{Q}_H, \qquad \forall \,\mathbf{q}_H\,\in\,\mathcal{Q}_H~.
\end{equation}
This assumption holds for the planar non-shearable rod that has been adopted as example, for which $Q_H$ is taken as the space of $C^1$-conforming piecewise cubic polynomials (for both components of $\mathbf{q}=\mathbf{m}$). 

We also assume that, for each configuration $\mathbf{q}_H\in Q_H$ there exist finite element spaces $V_h(\mathbf{q}_H)\subset X_{\mathbf{q}_H}$ and $\Pi_h(\mathbf{q}_H)\subset \Pi_{\mathbf{q}_H}$. The parameters $H$ and $h$ represent the sizes of the meshes of the rod (one-dimensional) and of the fluid ($d$-dimensional), respectively. The construction of $V_h(\mathbf{q}_H)$ and $\Pi_h(\mathbf{q}_H)$ requires a mesh over $\Omega_f(\mathbf{q}_H)$. Notice that meshes corresponding to different times, which will in general have different $\mathbf{q}_H$'s, are totally independent. To each $ \mathbf{r}_H\,\in\,Q_H$ we associate the fields $\mathbf{u}_h=\mathcal{U}_{h,\mathbf{q}_H,\mathbf{r}_H}$ and $p_h=\mathcal{P}_{h,\mathbf{q}_H,\mathbf{r}_H}$ unique solutions of
  \begin{equation}
    \left \{
    \begin{array}{l} (\mathbf{u}_h,p_h)\,\in\,V_h(\mathbf{q}_H)\times \Pi_h(\mathbf{q}_H),\\ \\
      \mathbf{u}_h = \mathcal{I}_h H_{\mathbf{q}_H} \mathbf{r}_H\qquad \mbox{ on }\mathcal{S}(\mathbf{q}_H) \\ \\
      \int_{\Omega_f(\mathbf{q}_H)} 2\mu \boldsymbol{\varepsilon}(\mathbf{u}_h)\bm{:} \boldsymbol{\varepsilon}({\mathbf w}_h)~d\Omega - \int_{\Omega_f(\mathbf{q}_H)} p_h \,\nabla \bm{\cdot}\,\mathbf{w}_h ~d\Omega
      = 0 \\ \\
      \int_{\Omega_f(\mathbf{q}_H)} m_h\,\nabla \bm{\cdot} \mathbf{u}_h ~d\Omega + \int_{\Omega_f(\mathbf{q}_H)} \tau_h \nabla p_h \cdot \nabla m_h~d\Omega = 0~
      \end{array}
    \right .
    \label{eq41}
    \end{equation}
 for all $m_h\,\in\,\Pi_h({\mathbf{q}_H})$ and all $\mathbf{w}_h$ in
\begin{equation}
V_{h0}({\mathbf{q}_H}) = \left \{ \mathbf{w}_h\,\in\,V_h(\mathbf{q}_H) ~|~ \mathbf{w}_h=0 \mbox{ on }\mathcal{S}(\mathbf{q}_H) \right \}~.
\end{equation}
Notice that if $\mu$ depends on the local strain rate then problem (\ref{eq41}) is nonlinear. The term containing $\tau_h$ is a stabilization term so that equal-order elements are rendered convergent by adopting $\tau_h = h^2/4\mu$. The operator $\mathcal{I}_h$ interpolates a vector function defined on $\mathcal{S}(\mathbf{q}_H)$ onto the trace space of $V_h(\mathbf{q}_H)$. We assume that the spaces $V_h$ and $\Pi_h$ are built from the Taylor-Hood $P_2-P_1$ element (in which case $\tau_h = 0$) or from the equal order $P_1-P_1$ element, so that it holds that (\cite{Taylor1973, Franca1992})
\begin{equation}
  \|\mathbf{u}_h - \mathcal{U}_{\mathbf{q}_H,\mathbf{r}_H}\|_{H^1(\Omega_f)}
  + \|p_h - \mathcal{P}_{\mathbf{q}_H,\mathbf{r}_H}\|_{L^2(\Omega_f)}
  \leq c\,h^k~,
  \label{eqerrstokes}
\end{equation}
with $c$ not depending on $h$ and $k=2$ for the $P_2-P_1$ case and $k=1$ for the $P_1-P_1$ one.

In our implementation we adopt Lagrange finite elements, so that $\mathcal{I}_h$ is the nodal interpolation at the mesh nodes that lie on $\mathcal{S}(\mathbf{q}_H)$. The extension operator $L_{\mathbf{q}_H}$ is the simplest possible, namely all nodal values at nodes not belonging to $\mathcal{S}(\mathbf{q}_H)$ are set to zero. 
We define a vector field $\mathbf{w}_h=\widehat{H}_{h,\mathbf{q}_H}\delta \mathbf{q}_H$ as taking the values $H_{\mathbf{q}_H}\delta \mathbf{q}_H$ at nodes that lie on $\mathcal{S}(\mathbf{q}_H)$, and zero at all others, i.e.,
\begin{equation}
\mathbf{w}_h=\widehat{H}_{h,\mathbf{q}_H}\delta \mathbf{q}_H\,\in\,V_h,\qquad
\mathbf{w}_h = \left \{ \begin{array}{ll}
\mathcal{I}_hH_{\mathbf{q}_H}\delta \mathbf{q}_H & \mbox{ on $\mathcal{S}(\mathbf{q}_H)$, }\\
\mathbf{0} & \mbox{ at all interior nodes of $V_h$.}
\end{array} \right.
\end{equation}

With the previous definitions, and using superscript $n$ to denote the time step, we will consider the following numerical methods:

\medskip

\noindent{\bf Fully discrete formulations of problem P2:} ``Find $\mathbf{q}_H^n\,\in\,Q_H$, for $n=0,1,2,\ldots$, such that $\mathbf{q}_H^0=\mathcal{I}_H \mathbf{q}_0$ and
\begin{equation}
  b_h\left (\mathbf{q}_H^n;\frac{\mathbf{q}_H^{n+1}-\mathbf{q}_H^n}{\delta t},\delta \mathbf{q}_H\right )+
  \,\langle D_{\mathbf{q}}E(t_*,\mathbf{q}_*),\delta \mathbf{q}_H \rangle
  + \, S(\mathbf{q}_H^n;\mathbf{q}_H^{n+1}-\mathbf{q}_H^n,\delta \mathbf{q}_H)
  - \,\langle \mathbf{F}(t_*),\delta \mathbf{q}_H \rangle = 0
  \label{eq42}
\end{equation}
for all $\delta \mathbf{q}_H\,\in\,Q_H$, where $b_h(\mathbf{q}_H;\cdot,\cdot):Q_H\times Q_H \to \mathbb{R}$ is given by 
\begin{equation}
  b_h(\mathbf{q}_H;\mathbf{p}_H,\delta \mathbf{q}_H)=
  \int_{\Omega_f(\mathbf{q}_H)} \left [ 2\mu\,\boldsymbol{\varepsilon} (\mathcal{U}_{h,\mathbf{q}_H,\mathbf{p}_H}) \bm{:} \boldsymbol{\varepsilon} (\widehat{H}_{h,\mathbf{q}_H}\delta \mathbf{q}_H)
-\mathcal{P}_{h,\mathbf{q}_H,\mathbf{p}_H}\,\nabla \cdot (\widehat{H}_{h,\mathbf{q}_H}\delta \mathbf{q}_H)
   \right ] ~d\Omega~,
\end{equation}
$\delta t = t_{n+1}-t_n$ is the time step and $S(\mathbf{q}_H;\cdot,\cdot):Q_H\times Q_H \to \mathbb{R}$ is a stabilizing bilinear form which will in general depend on the specific rod model.

The following choices for $t_*$ and $q_*$ define the different variants considered in our study:
\begin{itemize}
\item {\bf Explicit method:} $t_*=t_n$, $\mathbf{q}_*=\mathbf{q}_H^n$, $S=0$.
\item {\bf Semi-implicit method:} $t_*=t_n$, $\mathbf{q}_*=\mathbf{q}_H^n$, $S= S^{\,n}_{\mathrm{material}}\neq 0$.
\item {\bf Pseudo-implicit method:} $t_*=t_{n+1}$, $\mathbf{q}_*=\mathbf{q}_H^{n+1}$, $S=0$.''
\end{itemize}

Notice that in all three methods the bilinear (or nonlinear) form $b_h$ is evaluated at $\mathbf{q}_H^n$, which means that the fluid's geometry and mesh correspond to the previous time, $t_n$. The explicit method is a simple forward Euler of problem P2. The pseudo-implicit method, on the other hand, considers the elasticity of the rod {\bf implicitly} and thus requires nonlinear iterations to compute $\mathbf{q}_H^{n+1}$. It is however not an implicit method since, as just said, the fluid's treatment is explicit. Finally the semi-implicit method leads, when the fluid is Newtonian, to a linear algebraic system of equations, just as the explicit method. However, increased temporal stability is achieved by selecting $S$ equal to $S^{\,n}_{\mathrm{material}}$ (the material stiffness, which depends on the rod model) evaluated at configuration $\mathbf{q}_H^n$, which is a positive semidefinite bilinear form.

Assuming $\mathbf{q}_H^n$ known, the left-hand side of (\ref{eq42}) can be viewed as a residual $\mathcal{R}:\mathbb{R}^N\to \mathbb{R}^N$, $N$ being the dimension of $Q_H$, i.e.,
\begin{equation}
  \mathcal{R}(\mathbf{q}_H^{n+1}) = 0~. \label{eq44}
\end{equation}
The evaluation of the residual $\mathcal{R}$ is crucial for any iterative strategy to solve (\ref{eq44}). Each component of $\mathcal{R}$ corresponds to taking $\delta \mathbf{q}_H$ equal to one of the basis functions of $Q_H$. The computation is quite standard for all the terms in $\mathcal{R}$ except for the first one, which materializes the fluid-rod interaction through the bilinear (or nonlinear) form $b_h(\mathbf{q}_H^n;\cdot,\cdot)$.

For each candidate solution $\mathbf{q}_H^*$, the computation of $b_h(\mathbf{q}_H^n;(\mathbf{q}_H^{*}-\mathbf{q}_H^n)/\delta t,\delta \mathbf{q}_H)$ involves the following steps:
\begin{enumerate}
\item Build a mesh $\mathcal{T}_h$ for the fluid domain $\Omega_{f}^n=\Omega_f(\mathbf{q}_H^n)$ with the geometry defined by $\mathbf{q}_H^n$.
  \item Solve the fluid problem (\ref{eq41}) with $\mathbf{r}_H=(\mathbf{q}_H^{*}-\mathbf{q}_H^n)/\delta t$. Let $(\mathbf{u}_h^*,{p}_h^*)$ be the corresponding velocity-pressure pair. 
  \item For each $i=1,\ldots,n_{\mbox{\small{dof}}}$, if $\mathbf{N}^i$ is the $i-$th basis function of $Q_H$, compute $H_{\mathbf{q}_H^n}\mathbf{N}^i$ at the {\bf fluid nodes} on the rod's boundary. Setting it to zero for all the other fluid nodes, one obtains $\mathbf{w}_h^i=\widehat{H}_{h,\mathbf{q}_H^n}\mathbf{N}^i$.
  \item The first term in the $i-$th component of $\mathcal{R}(\mathbf{q}_H^*)$ is then obtained from the integral
    \begin{equation}
      \mathcal{R}^i(\mathbf{q}_H^*) = b_h\left ( \mathbf{q}_H^n;\frac{\mathbf{q}_H^{*}-\mathbf{q}_H^n}{\delta t},\mathbf{N}^i\right ) + \mbox{o.t.} =
            \int_{\Omega_f(\mathbf{q}_H^n)} \left [ 2\mu\,\boldsymbol{\varepsilon} (\mathbf{u}_h^*) \bm{:} \boldsymbol{\varepsilon} (\mathbf{w}_h^i)
-{p}_h^*\,\nabla \cdot \mathbf{w}_h^i
   \right ] ~d\Omega~ + \mbox{o.t.}, \label{eq:residual}
      \end{equation}
where ``o.t.'' stands for the ``other terms'' in the left-hand side of (\ref{eq42}), which all come from the solid model. 
\end{enumerate}

Notice that steps 1 and 3 do not depend on $\mathbf{q}_H^*$ and can be computed just once per time step. It remains to provide the discretized equations of the rod so as to compute the remaining terms in $\mathcal{R}^i(\mathbf{q}_H^*)$. These equations, which depend on the rod model, are provided in detail for the case of the planar non-shearable rod in the following section.

%\medskip

\noindent{\bf Remark:} For any $\mathbf{q}\in\mathcal{Q}$, any $\mathbf{r}_H\in Q_H$ and any basis function $\mathbf{N}^i\in Q_H$, the difference between $b(\mathbf{q};\mathbf{r}_H,\mathbf{N}^i)$ and $b_h(\mathbf{q};\mathbf{r}_H,\mathbf{N}^i)$ is that in the latter the exact Stokes solution $(\mathcal{U}_{\mathbf{q},\mathbf{r}_H},\mathcal{P}_{\mathbf{q},\mathbf{r}_H})$ is approximated on a mesh of $\Omega_f(\mathbf{q})$ of size $h$. A first formal estimate of the consistency error is
\begin{equation*}
b(\mathbf{q};\mathbf{r}_H,\mathbf{N}^i)-
b_h(\mathbf{q};\mathbf{r}_H,\mathbf{N}^i)=
\end{equation*}
\begin{equation*}
=\int_{\Omega_f(\mathbf{q})} \left [ 2\mu\,\boldsymbol{\varepsilon} (\mathcal{U}_{\mathbf{q},\mathbf{r}_H}) \bm{:} \boldsymbol{\varepsilon} (\widehat{H}_{\mathbf{q}} \mathbf{N}^i)
- \mathcal{P}_{\mathbf{q},\mathbf{r}_H}\,\nabla \cdot (\widehat{H}_{\mathbf{q}} \mathbf{N}^i)
   \right ] ~d\Omega~ -
\end{equation*}
\begin{equation*}
-\int_{\Omega_f(\mathbf{q})} \left [ 2\mu\,\boldsymbol{\varepsilon} (\mathcal{U}_{h,\mathbf{q},\mathbf{r}_H}) \bm{:} \boldsymbol{\varepsilon} (\widehat{H}_{h,\mathbf{q}} \mathbf{N}^i)
-\mathcal{P}_{h,\mathbf{q},\mathbf{r}_H}\,\nabla \cdot (\widehat{H}_{h,\mathbf{q}} \mathbf{N}^i)
   \right ] ~d\Omega~ \leq
\end{equation*}
\begin{equation*}
\leq c_1 \left (
  \|\mathcal{U}_{\mathbf{q},\mathbf{r}_H} - \mathcal{U}_{h,\mathbf{q},\mathbf{r}_H}\|_{H^1(\Omega_f)}
  + \|\mathcal{P}_{\mathbf{q},\mathbf{r}_H} - \mathcal{P}_{h,\mathbf{q},\mathbf{r}_H}\|_{L^2(\Omega_f)} \right )+ c_2 \left \|H_{\mathbf{q}}\mathbf{N}^i-\mathcal{I}_h H_{\mathbf{q}}\mathbf{N}^i)\right \|_{H^{1/2}(\mathcal{S})}  \leq
  \end{equation*}
  \begin{equation*}
  \leq c_1\,h^k + c_2'\,h^{3/2}~,
\end{equation*}
  where $c_1$, $c_2$ and $c_2'$ depend on $\mathbf{q}$, $\mathbf{r}_H$ and $\mathbf{N}^i$, but not on $h$, and we have used that $\mathbf{N}^i$ is in $H^2(0,S)$ in formally estimating the third term (interpolation error). By applying a duality argument the term $h^k$ can probably be improved to $h^{2k}$. Recall that the total error involves the consistency error just discussed together with the interpolation error coming from the (Hermite) solid model, $
  \|\mathbf{q}-\mathcal{I}_H\mathbf{q}\|_{H^2(0,S)} \leq c_3\,H^2
  $, where $\mathcal{I}_H\mathbf{q}\in V_H$ is the Hermite interpolant of $\mathbf{q}$.
  
%---------------------------------------------------------------------
%---------------------------------------------------------------------
  
\subsection{Software: Firedrake implementation for a planar non-shearable rod}

Let us discuss the ingredients involved in the finite element implementation
of the fluid structure interaction problem for the particular case of a planar non-shearable rod. 
The finite element library Firedrake \cite{Rathgeber2016} has been adopted.
This platform allows expressive specification of PDE's using the Unified Form Language (UFL)
from the FEniCS Project \cite{Alnaes2014} through a python interface.

\subsubsection{Kinematics, mesh generation and remeshing}

\begin{figure} 
   \begin{center}
       \scalebox{0.925}{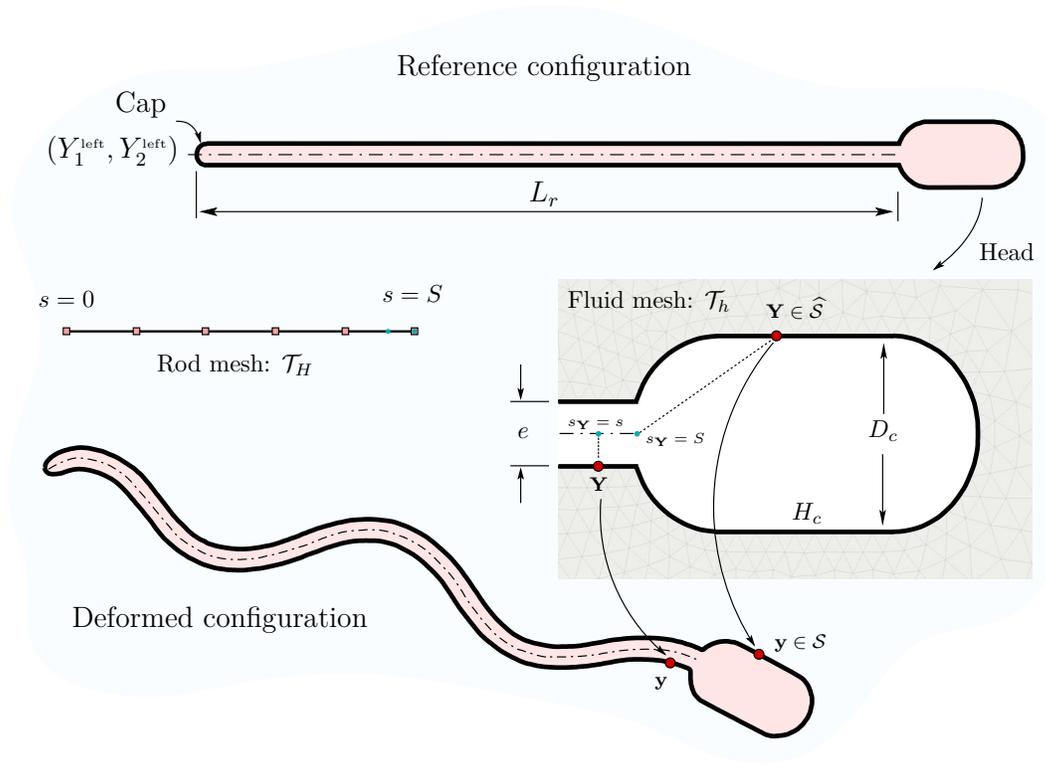}
       \caption{Schematic of the computational domain showing the reference and a deformed 
       configuration of the rod swimmer. The 2D fluid mesh and 1D rod mesh are also
       shown.}\label{fig:refswimmer}
   \end{center}
\end{figure}

For the implementation we first need to collect some results previously introduced
related to the kinematics of the planar non-shearable rod. 
Let us recall that the rod configuration at time $t_n$ is encoded in the $\mathbf{q}_H^n$
function, so that the fluid domain $\Omega_f(\mathbf{q}_H^n)$ can be built from $\mathbf{q}_H^n$ and the parameters $\{\theta_1, \theta_2\}$ computed at the reference configuration.
To fix ideas, let us consider the rod shown in Fig. \ref{fig:refswimmer}
and suppose we have constructed an initial mesh around it. This will be the 
reference configuration. For each node over the wet boundary $\widehat{\mathcal{S}}$
(the red dots in the figure) we record its position $\mathbf{Y}$ in this reference state.
From this, the parameter $s_\mathbf{Y} = \sigma(\mathbf{Y}) \in [0,S]$ 
to which each material point is attached can be found.
Along its evolution, in any configuration of the rod, the actual position $\mathbf{y}$
of the same material node is computed according to
\begin{equation}
 \mathbf{y} = \mathbf{q}_H(s_{\mathbf{Y}}) + \theta_1(\mathbf{Y}) \, \textbf{q}_H'(s_{\mathbf{Y}}) +
 \theta_2(\mathbf{Y}) \, a(s_{\mathbf{Y}}) \, \underbrace{\mathbf{J} \, \frac{\textbf{q}_H'(s_{\mathbf{Y}})}{\parallel \textbf{q}_H'(s_{\mathbf{Y}}) \parallel}}_{\mathbf{\check{n}}_H} \label{eq:updateWet}
\end{equation}
In the numerical experiments to be shown later on we consider two particular cases
\begin{equation}
a(s) = 1, ~~~~\mbox{and}~~~ a(s) = \frac{1}{\parallel \textbf{q}_H'(s) \parallel} \label{eq:aofs}
\end{equation}
Consider the rod shown in Fig. \ref{fig:refswimmer} whose length in the reference configuration
is denoted by $L_r$ (\textcode|Lref| in the code). This rod has a cap with rounded shape at the left end, 
a head at the right end and thickness $e \ll L_r$. The left cap is described by half of a 
circle of diameter $e$. The right head has two circular caps of diameter $D_c = \beta_d \, e$
and a central region of height $H_c = \beta_h \, e$, being $\beta_d,\beta_h > 1$ 
user defined parameters such that its size can be easily altered. 
These geometrical features are all indicated in the figure.
For illustration purposes the left and right parts are dealt differently in this example, 
namely, the left rounded tip is part of the rod's domain, i.e., all the material points
have $0 \le s_{\mathbf{Y}} < S$, while the material points belonging to the head are all attached to 
$s_{\mathbf{Y}} = S$. This amounts to define
\begin{equation}
 s_{\mathbf{Y}} = 
 \left \{ 
\begin{array}{ll}
\dfrac{S}{L_r}(Y_1 - Y_1^{\mbox{\tiny left}}) & ~\mbox{if}~ Y_1 - Y_1^{\mbox{\tiny left}} \in [0,S) \\
 \\
 S & ~\mbox{if}~ Y_1 - Y_1^{\mbox{\tiny left}} \ge S
\end{array}
 \right. \label{eq:sigmapart}
\end{equation}
where $( Y_1^{\mbox{\tiny left}}, Y_2^{\mbox{\tiny left}} )$ 
is the initial position of the left end.
Since $\mathbb{D} = \mathbf{I}$ and 
$\mathbf{M}(s) = (s\,L_r/S + Y_1^{\mbox{\tiny left}}, Y_2^{\mbox{\tiny left}})$
by Eq. (\ref{eq:titasgen}) the 
parameters $\theta_i$ are thus given by
\begin{equation}
(\theta_1(\mathbf{Y}),\theta_2(\mathbf{Y})) = 
 \left \{ 
\begin{array}{ll}
 (0,Y_2-Y_2^{\mbox{\tiny left}}) & ~\mbox{if}~s_\mathbf{Y} \in [0,S) \\
 \\
 (~Y_1-Y_1^{\mbox{\tiny left}} - L_r, ~Y_2-Y_2^{\mbox{\tiny left}}) & ~\mbox{otherwise}
\end{array}
 \right. \label{eq:titas}
\end{equation}
To update the position of nodes over the wet boundary, Eq. (\ref{eq:updateWet}) is implemented by the
python function \textcode|UpdateWetBoundary()|. This in turn calls some user defined functions, namely,
\textcode|SigmaofY()| (Eq. \ref{eq:sigmapart}), \textcode|ComputeThetas()| (Eq. (\ref{eq:titas})) 
and \textcode|aofs()| (Eq. (\ref{eq:aofs})) shown next. Additionally, the function \textcode|aofsp()| 
that computes $a'(s)$, to be used later on, is also shown.
\begin{lstlisting}[language=Python]
# Program these functions according to your problem
def SigmaofY(xnodref): # Eq. (53)
    sY = S*(xnodref[0] - Y1left)/Lref
    if(sY > S):
        sY = S
    return sY
def ComputeThetas(xnodref): # Eq. (54)
    sY = SigmaofY(xnodref)
    if(sY < S):
        theta1 = 0
    else:
        theta1 = xnodref[0] - Y1left - Lref
    theta2 = xnodref[1] - Y2left
    return theta1, theta2

def aofs(modqp): # Eq. (52)
    return 1/modqp
def aofsp(qp,qpp): # Derivative of (52)
    return -np.dot(qp,qpp) / np.linalg.norm(qp)**3

def UpdateWetBoundary(xref, qH): # Eq. (51)
    new_y = np.zeros(shape=(nnod_wet,2))
    for j in range(nnod_wet): # Swept over wet nodes
        sY = SigmaofY(xref[j,:])
        el_inside = int(one_minus_eps * sY / h_r)
        xi = 2*sY/h_r - (2*el_inside+1)
        q = EvalqH(qH, el_inside, xi)
        qp = EvalqHp(qH, el_inside, xi)
        theta1,theta2 = ComputeThetas(xref[j,:])
        modqp = np.linalg.norm(qp)
        ny = q + theta1*qp + theta2*aofs(modqp)/modqp*np.dot(Jota,qp)
        new_y[j,:] = ny[:]
    return new_y
\end{lstlisting}
Functions \textcode|EvalqH()| and \textcode|EvalqHp()| that specify the rod configuration 
are explained later on when discussing the discretization of the rod problem.
Table \ref{tab:symbols_python} lists several important objects and 
variables that appear in the code above and others to be presented afterwards.
Notice that some of these variables may have global scope in the code.

\begin{table}[ht!!]
\begin{center}
\caption{List of important symbols and variables that appear in the Firedrake-python implementation 
made available alongside the article.} \label{tab:symbols_python}
{\small 
\begin{tabularx}{\textwidth}{r|c|X}
{\bf Symbol} & {\bf Python}  &  {\bf Description} \\
\hline
$L_r$               & \textcode|Lref|     &  Length of rod in the reference configuration $(\mathbf{M}, \mathbb{D})$ \\
$( Y_1^{\mbox{\tiny left}}, Y_2^{\mbox{\tiny left}} )$ & \textcode|Y1left, Y2left|   & Position coordinates of the left end in the reference configuration \\
$(C_{\epsilon},C_{\kappa})$  & \textcode|Ce,Ck|   & Rod's cross-sectional elastic properties\\
$\mathbf{I}$        & \textcode|iden, I|     & \textcode|np.array([[1,0],[0,1]]), Identity(2)| - Identity matrix \\
$\mathbf{J}$        & \textcode|Jota, J|     & \textcode|np.array([[0,-1],[1,0]]), as\_matrix([[0, -1], [1, 0]])| - Unitary skew--symmetric matrix \\
$\mathcal{T}_H$     & \textcode|mesh\_r|  & Firedrake object - 1D rod mesh \\
$N_{\mbox{\tiny el}}$ & \textcode|Nelrod| & Number of elements in rod mesh \\
$N_{\mbox{\tiny nod}}$ & \textcode|Nnodrod|& \textcode|=Nelrod+1|, Number of nodes in rod mesh \\
$N_{\mbox{\tiny dofs}}$ & \textcode|Ndofs|& \textcode|=4*Nnodrod|, Total number of rod's dofs  \\
$Q_H$               & \textcode|QH|     & \textcode|VectorFunctionSpace(mesh\_r, 'HER', degree=3, dim=2)| - Hermite vector function space for rod discretization  \\
$\mathbf{N}^i$      & \textcode|deltaqH[i-1]|  & \textcode|=Function(QH)| - Nodal basis functions for the space $Q_H$ \\
$H$                 & \textcode|h\_r|     & \textcode|=Lref/Nelrod|, Fixed elemental length in the rod mesh \\
$s_{\mathbf{Y}} = \sigma(\mathbf{Y})$     & \textcode|sY| & \textcode|=SigmaofY(xnodref)| - Rod coordinate to which point $\mathbf{Y}$ is attached \\
$\mathbf{q}^{n+1}_H$    & \textcode|qH|     & \textcode|=Function(QH)| - Rod configuration at current time as function belonging to the space $Q_H$ \\
$\mathbf{q}^n_H$        & \textcode|qHn|     & \textcode|=Function(QH)| - Rod configuration at previous time as function belonging to the space $Q_H$ \\
${\mathbf{\dot{q}}}^{n+1}_H$    & \textcode|qHdot|     & \textcode|=Function(QH)| - Rod velocity at current time as function belonging to the space $Q_H$ \\
$\underline{\alpha}=\{\alpha_i\}_{i=1}^{N_{\mbox{\tiny dofs}}}$    & \textcode|alpha|     & \textcode|=np.zeros(Ndofs)| - Rod velocity dofs at current time\\
$\mathbf{q}_H(s)$      & \textcode|q|       & Rod configuration evaluated at arbitrary $s$ \\
$\mathbf{q}_H'(s)$     & \textcode|qp|      & First derivative of $\mathbf{q}_H$ evaluated at arbitrary $s$\\  
$\mathbf{q}_H''(s)$    & \textcode|qpp|     & Second derivative of $\mathbf{q}_H$ evaluated at arbitrary $s$\\
$\delta{t}$    & \textcode|dt|     & Time step \\
$-$    & \textcode|EM| & \textcode|=MixedElement([VectorElement('Lagrange','triangle',1),| \textcode|  FiniteElement('Lagrange','triangle',1)])|   \\
$V_{h0}\times \Pi_h$  & \textcode|W\_f|       & \textcode|=FunctionSpace(mesh\_f, EM)| - Mixed Lagrangian function space for fluid discretization \\
$\mathcal{T}_h$       & \textcode|mesh\_f|    & Firedrake object - 2D fluid mesh \\
$\mathbf{w}_h^i$      & \textcode|wh[i-1]|    & \textcode|=Function(W\_f.sub(0))|, $i$--th induced velocity basis function \\
$(\mathring{\mathbf{u}}_h,p_h)$               & \textcode|uo,p|      & \textcode|=split(sol\_f)|, \textcode|sol\_f = Function(W\_f)| - velocity--pressure pair \\
${\mathbf{u}}_h$      & \textcode|u|           & Final velocity - ${\mathbf{u}}_h = \mathring{\mathbf{u}}_h + \sum_{i=1}^{N_{\mbox{\tiny dof}}} { \alpha_i^{n+1} \, \mathbf{w}_h^i }$ \\
$(\mathbf{v}_h,q_h)$  & \textcode|v,q|         & \textcode|=TestFunctions(W\_f)| \\
$\tau_h$          & \textcode|tauh|            & \textcode|=hk*hk/(4*mu)| (\textcode|hk = CellDiameter(mesh\_f)|) - Stabilization parameter in equal-order formulation \\
$-$               & \textcode|refinement|      & Typical element size close to the wet boundary \\
$-$               & \textcode|nnod\_wet|       & Number of nodes over the wet boundary \\
$-$               & \textcode|scaling|         & \textcode|=np.array([1.,h\_r])| - Scaling factors of Hermite shape functions \\
$-$               & \textcode|one\_minus\_eps| & \textcode|1-1e-12|, geometrical tolerance \\
\end{tabularx}
}
\end{center}
\end{table}
The user defined functions can be adapted to consider other cases. 
For instance, for the simpler case of a rod with both 
flats ends the parameters $\theta_i$ are just given by
$$
(\theta_1(\mathbf{Y}),\theta_2(\mathbf{Y})) = (0, ~Y_2-Y_2^{\mbox{\tiny left}})~~\forall s_{\mathbf{Y}} \in [0,S]
$$
For the sake of simplicity we have ommited additional arguments, such as $\mathbb{D}$ and $\mathbf{M}$,
that can be passed to these functions if necessary to deal with different or more complex
reference configurations.

\medskip

Once the new position of nodes over the wet boundary has been found, a finite element
mesh is built. Linear triangles are used for discretization of the fluid domain.
The mesh generator gmsh \cite{Geuzaine2009} (version 3.0.6) through the python 
library pygmsh \cite{pygmsh2020} (version 6.1.0) is used for such a task. 
The advancing front algorithm from the gmsh suite was chosen in the numerical
experiments below (\textcode|gmsh_args=['-algo','front2d']|).
Two functions are used to handle the mesh generation in our code, namely,
\textcode|MeshFluidRegion()| and \textcode|RemeshFluidRegion()|.
The former can be adapted so as to define different initial 
rod configurations. 
\begin{lstlisting}[language=Python]
...

mesh_f, ... = MeshFluidRegion( ..., refinement, gmsh_args, ... )
.
.
.
if(need_remesh):
  mesh_f, ... = RemeshFluidRegion( mesh_f, x_f, xref, . , list_ord, ... )

...
\end{lstlisting}
Many of the arguments in those functions have been ommited
to avoid non-essential details. They nevertheless appear in the source 
code made available alongside the article.
The characteristic size $h$ of elements close to the rod's boundary $\mathcal{S}(\mathbf{q}_H^0)$
in the initial configuration is controlled by the parameter \verb+refinement+. 
A key issue during remeshing is to preserve the position
of nodes located over the wet boundary. To that end,
the mesh coordinates \verb+x_f+ of the original mesh are passed as an additional
argument to the remeshing function along an ordered list of nodes \verb+list_ord+ to
be connected during remeshing to form the new wet boundary. 
Since we expect to use relatively large time steps, remeshing will be applied at
each time in the numerical experiments to be shown later on.
Fig. \ref{fig:refswimmer} shows a detail of a typical initial mesh used in computations.

\subsubsection{Discrete rod problem}

Consider the rod domain $[0,S]$ discretized into $N_{\mbox{\tiny el}}$ elements of uniform length
$H = S / N_{\mbox{\tiny el}}$. These parameters are identified in the implementation
by \textcode|Nelrod| and \textcode|h_r|, respectively.
Define the unknown discrete rod velocity $\mathbf{\dot{q}}_H \in Q_H$ 
between $t_n$ and $t_{n+1}$ as
\begin{equation}
 \mathbf{\dot{q}}_H^{n+1} = \frac{\mathbf{q}_H^{n+1} - \mathbf{q}_H^{n}}{\delta{t}}
\end{equation}
which is writen as a linear combination of the nodal basis functions 
$\left \{ {\bf N}^i \right \}_{i=1}^{N_{\mbox{\tiny dof}}}$ of $Q_H$, i.e.,
\begin{equation}
 \mathbf{\dot{q}}_H^{n+1} = \sum_{i=1}^{N_{\mbox{\tiny dof}}} { \alpha_i \, {\bf N}^i }
\end{equation}
For discretization of the rod we use, as anticipated, one dimensional $C^1$-conforming Hermite elements and
a regular grid with elements of length $H$.
Being $\mathbf{\dot{q}}_H^{n+1}$ a vector field with two components in 2D, the total number of dofs in the rod is $N_{\mbox{\tiny dof}} = 4\,N_{\mbox{\tiny nod}}$ with $N_{\mbox{\tiny nod}} = N_{\mbox{\tiny el}} + 1$.
On each element $e=1, \dots, N_{\mbox{\tiny el}}$ of $\mathcal{T}_H$ we approximate any function of $Q_H$ as 
\begin{equation}
    \mathbf{q}_H(\len)|_e = \underline{\nmat}^e(\xi_s) \, \overline{\mathbf{q}}^{\,e}_H,~~s = s_e \, \frac{1-\xi_s}{2} + s_{e+1} \,\frac{1+\xi_s}{2},~s_e = e \, H
\end{equation}
with the nodal values being $\overline{\mathbf{q}}^{\,e}_H$ and 
the matrix of shape functions defined by 
\begin{equation}
    \underline{\nmat}^{e}(\xi) =
    \begin{pmatrix}
    \eye N_{\posname}^{-}(\xi) & \frac{H}{2} \eye N_{\tanname}^{-}(\xi) & \eye N_{\posname}^{+}(\xi) & \frac{H}{2} \eye N_{\tanname}^{+}(\xi)
    \end{pmatrix} ,
\end{equation}
with the standard Hermite's polynomials
\begin{subequations}
    \begin{align}
    N_{\posname}^{-}(\xi) & = \frac14 (2+\xi)(1-\xi)^2 ,\\
    N_{\tanname}^{-}(\xi) & = \frac14 (1+\xi)(1-\xi)^2 ,\\
    N_{\posname}^{+}(\xi) & = \frac14 (2-\xi)(1+\xi)^2 ,\\
    N_{\tanname}^{+}(\xi) & = -\frac14 (1-\xi)(1+\xi)^2 .
    \end{align}
\end{subequations}
These shape functions are defined for $\xi\in[-1,+1]$.
Their evaluation at arbitrary points, as done for instance
in \textcode|EvalqH()|, needs to be programmed manually since no Firedrake function
that works for Hermite elements is available for such a task.

\medskip

The following code block shows the definition of the vector space 
$Q_H$ as a Firedrake object. Also, in our implementation we explicitly build 
the basis functions ${\bf N}^i$ of the space $Q_H$. These functions are 
used later on as test functions in (\ref{eq:residual}). In the code
we also illustrate how any initial condition (a planar rod of lenght $L_r$)
is defined by using the Firedrake function \textcode|project()|.
\begin{lstlisting}[language=Python]
mesh_r = IntervalMesh(Nelrod, S)
QH = VectorFunctionSpace(mesh_r, 'HER', degree=3, dim=2)
count, deltaqH = 0, []
for i in range(Nnodrod):
    [deltaqH.append(Function(QH)) for d in range(4)]
    deltaqH[count+0].dat.data[2*i][0]   = 1
    deltaqH[count+1].dat.data[2*i][1]   = 1
    deltaqH[count+2].dat.data[2*i+1][0] = h_r
    deltaqH[count+3].dat.data[2*i+1][1] = h_r
    count = count + 4
# Initial condition
x_r = SpatialCoordinate(mesh_r)
qH  = project(as_vector([x_r[0]*Lref/S + Y1left, Y2left]),QH)
\end{lstlisting}

\medskip

Let us now see how the terms arising from the rod model in Eq. (\ref{eq42}) are computed for the case of a planar non-shearable rod. So, the terms we are considering here are
\begin{equation}
  \mbox{Terms from rod model }=\,\langle D_{\mathbf{q}}E(t_*,\mathbf{q}_*),\delta \mathbf{q}_H \rangle
  + \, S(\mathbf{q}_H^n;\mathbf{q}_H^{n+1}-\mathbf{q}_H^n,\delta \mathbf{q}_H)
  - \,\langle \mathbf{F}(t_*),\delta \mathbf{q}_H \rangle~,
  \label{eq:terms}
\end{equation}
where $t_*=t_n$, $\mathbf{q}_*=\mathbf{q}_H^n$ correspond to the explicit method, for which $S=0$, and to the semi-implicit method, for which $S$ is given by
\begin{eqnarray}
S(\mathbf{q}_H^n; \mathbf{{q}}_H^{n+1}-\mathbf{{q}}_H^{n},\delta \mathbf{q}_H)
&=&\int_0^S \delta \mathbf{q}_H\cdot \theoperator{K}_{\mathrm{material}}^{\, n}\,(\mathbf{{q}}_H^{n+1}-\mathbf{{q}}_H^{n})~\difflen \nonumber \\
&=&
 \delta{t} \, \int_0^S{\left ( \theoperator{B}^{\,n} \delta{\mathbf{q}_H} \right ) \cdot \mathbf{C} \, \left ( \theoperator{B}^{\,n} \mathbf{\dot{q}}_H^{n+1} \right) }~\difflen~\nonumber \\
 & = & \delta{t}\,\int_0^S \delta \elasstrain \cdot \mathbf{C}\,\Delta \elasstrain~\difflen.
 \end{eqnarray}
 Above, $\delta \elasstrain= \theoperator{B}^{\,n} \delta{\mathbf{q}_H}$ and $\Delta \elasstrain = \theoperator{B}^{\,n} \mathbf{\dot{q}}_H^{n+1}$.
 Notice that in the notation we have used a supraindex $n$ to denote the operators $\theoperator{B}$ and $\theoperator{K}_{\mathrm{material}}$ evaluated at $({\mathbf{q}_H^{n}}', {\mathbf{q}_H^{n}}'')$. Also notice for later use that
 \begin{equation}
 S(\mathbf{q}_H^n; \mathbf{{q}}_H^{n+1}-\mathbf{{q}}_H^{n},\delta \mathbf{q}_H)=
\delta t\, S(\mathbf{q}_H^n; \mathbf{\dot{q}}_H^{n+1},\delta \mathbf{q}_H)
 \end{equation}
 In turn, the pseudo-implicit method considers $t_*=t_{n+1}$, $\mathbf{q}_*=\mathbf{q}_H^{n+1}$, and $S=0$. The first term in (\ref{eq:terms}) can also be written down explicitly, i.e.,
 \begin{equation}
 \langle D_{\mathbf{q}}E(t_*,\mathbf{q}_*),\delta \mathbf{q}_H \rangle
 =~\int_0^S{\left ( \theoperator{B}^{\,n} \delta{\mathbf{q}_H} \right) \cdot \mathbf{C} \, \elasstrain^{\,*} }~\difflen~=~\int_0^S \delta \elasstrain \cdot \mathbf{C} \, \elasstrain^{\,*}~\difflen,
 \end{equation}
 where $\elasstrain^{\,*}$ is the elastic strain corresponding to time $t_*$ and  configuration $\mathbf{q}_*$.

The implementation shown below is flexible and allows one to choose among the three different methods just proposed. Care was taken that there is
a close correspondence between the presented formulation and the code of its 
Firedrake implementation. 
Functions \textcode|epsilon()| and \textcode|kappa()| 
compute the strain measures given in (\ref{eq:epsilon}) and (\ref{eq:kappa}),
and \textcode|depsilon()| and \textcode|dkappa()| their
first variations with respect to \textcode|dq| given
by Eqs. (\ref{eq:depsilon}) and (\ref{eq:dkappa}). 
Lastly, \textcode|b_RodTerms_b()|
computes all three terms of (\ref{eq:terms}). In particular, \textcode|qs| corresponds to $\mathbf{q}_*$, the variables \textcode|eps| and \textcode|kap| are the components of $\elasstrain^*$, the variables \textcode|deps| and \textcode|dkap| are the components of $\delta \elasstrain$ and the variables \textcode|depsd| and \textcode|dkapd| are those of $\Delta \elasstrain$.

Notice that the generalized force considered in the implementation is of the form
\begin{equation}
\langle \mathbf{F}(t), \delta \mathbf{q} \rangle = \int_0^S \left ( \bm{f}_0(s,t)\cdot \delta \mathbf{q}(s)+ \bm{f}_1(s,t)\cdot \delta \mathbf{q}'(s)\right )~\difflen
\end{equation}
where $\bf{f}_0$ and $\bf{f}_1$ are user defined functions.

\begin{lstlisting}[language=Python]
def epsilon(t,qp,L0):
    return (dot(t,qp) - L0)
def kappa(t,np,k0):    
    return (dot(t,np) - k0)
def depsilon(t,dq):
    return dot(t,dq.dx(0))
def dkappa(t,n,qp,qpp,dq):
    dnttmJ = 2*outer(n,t) - J
    dqp  = dq.dx(0)
    dqpp = dq.dx(0).dx(0)
    mqp2 = dot(qp,qp)
    aux  = dot(dnttmJ,dqp)
    return dot(qpp,aux)/mqp2 - dot(n,dqpp)/sqrt(mqp2)

def b_RodTerms_b(dt, eref, espont, fe, qHn, qHdot, dq):
    if(scheme == 'explicit' or scheme == 'semi_implicit'):
        qs = qHn
    else:
        qs = qHn + Constant(dt)*qHdot
    qp  = qs.dx(0)
    qpp = qs.dx(0).dx(0)
    mqp = sqrt(dot(qp,qp))  
    t   = qp/mqp 
    n   = dot(J,t)
    Bper = I - outer(t,t)
    tp   = dot(Bper,qpp)/mqp
    np   = dot(J,tp)
    deps = depsilon(t,dq)
    eps  = epsilon(t,qp,eref[0]+espont[0])
    dkap = dkappa(t,n,qp,qpp,dq)
    kap  = kappa(t,np,eref[1]+espont[1])
    Fgrod = deps*Ce*eps + dkap*Ck*kap
    if(scheme == 'semi_implicit'): # Stabilization
        depsd = depsilon(t,qHdot)
        dkapd = dkappa(t,n,qp,qpp,qHdot)
        Fgrod += Constant(dt)*(depsd*Ce*deps + dkapd*Ck*dkap)
    Fgrod -= dot(fe[0],dq) + dot(fe[1],dq.dx(0))
    return assemble(Fgrod * dx(degree=5))
\end{lstlisting}

\subsubsection{Induced velocity basis functions}

Construction of the velocity basis functions $\mathbf{w}_h^i$ deserves special attention
since it is a key ingredient in the proposed implementation and is related
to the construction of the space of kinematically admissible velocity fields
in the discrete case.
Consider the finite element partition $\mathcal{T}_h$ of $\Omega_f(\mathbf{q}_H^n)$.
The $i$--th basis function correponds to
taking $\mathbf{N}^i,~i=1,\dots,N_{\mbox{\tiny dof}}$ in Eq. (\ref{eq:inducedvel}).
For a given node over the wet surface $\mathcal{S}(\mathbf{q}_H^n)$, with coordinate
$\mathbf{y}_j$, corresponding to the material point $\mathbf{Y}$ in the reference configuration
and attached to the rod coordinate $s_{\mathbf{Y}}$, we have
\begin{equation}
\mathbf{w}_h^i(\mathbf{y}_j) = \left ({H}_{\mathbf{q}_H^n} \mathbf{N}^i\right )(\mathbf{y}_j) = \mathbf{N}^i(s_{\mathbf{Y}})+ \left [ \theta_1(\mathbf{Y})\, \mathbf{I} + \theta_2(\mathbf{Y})\,\mathbf{G}(s_{\mathbf{Y}}) \right ] \, [\mathbf{N}^i(s_{\mathbf{Y}})]'~. \label{eq:whcode}
\end{equation}
and is defined as zero at all the nodes of $\Omega_f(\mathbf{q}_H^n)$ that are not on the wet surface. These nodal values define
$\mathbf{w}_h^i$, which in the discrete formulation appears as $\mathbf{w}_h^i=\widehat{H}_{h,\mathbf{q}_H^n} \mathbf{N}^i$.
The revelant code to compute these functions is shown in the block below. 
For the $i$--th dof we compute the function values of \textcode|w_h[i-1]| 
defined as a function belonging to the velocity space \textcode|W_f.sub(0)|.
The python function \textcode|ApplyHhatqN()| executes a sweep over all rod's 
dofs, through four nested loops (over \textcode|i,d,c|), 
thus totalizing \textcode|2*2*Nnodrod| functions. The velocity is computed
at each node on the wet surface by sweeping over the \textcode|nnod_wet| nodal 
coordinates \textcode|xref[j,:]| and calling function \textcode|VelWet()|, 
which is also shown. The \textcode|SupportInfo()| function provides auxiliary 
information to verify if a given point \textcode|j| is in the support of
the corresponding basis function \textcode|i|.

\begin{lstlisting}[language=Python]
def VelWet(xnodref,qp,qpp,Ndot,Ndotp): # Eq. (61)
    modqp = np.linalg.norm(qp)
    tan  = qp / modqp
    nn = np.dot(Jota,tan)
    nntensorqp =  np.tensordot(nn,qp,0)
    Bper = iden - np.tensordot(tan,tan,0)
    theta1,theta2 = ComputeThetas(xnodref)
    vel_wetnod = Ndot + theta1*Ndotp + \
                 theta2*aofs(modqp)/modqp*np.dot(Jota*Bper,Ndotp) + \
                 theta2*aofsp(qp,qpp)/modqp*np.dot(nntensorqp,Ndotp)
    return vel_wetnod

def ApplyHhatqN(W_f, xref, list_nodes_wet, qH): 
    scaling = np.array([1.,h_r])
    Ndot, Ndotp = np.zeros(2), np.zeros(2)
    count, wh = 0, []
    for i in range(Nnodrod): # Swept over rod nodes
        [wh.append(Function(W_f.sub(0), name='wh')) for d in range(4)]
        for j in range(nnod_wet): # Swept over wet fluid nodes
            k = list_nodes_wet[j]
            el_inside, xi, idn = SupportInfo(i, xref[j,:])
            if(el_inside == i-1 or el_inside == i):
                qp  = EvalqHp(qH, el_inside, xi)
                qpp = EvalqHpp(qH, el_inside, xi)
                HermBase, DerHermBase = HermB(xi), DHermB(xi)
                for d in range(2): # Swept over dofs per rod node
                    for c in range(2):
                        Ndot.fill(0.), Ndotp.fill(0.)
                        Ndot[c]  = scaling[d]*HermBase[idn[d]]
                        Ndotp[c] = scaling[d]*(2/h_r)*DerHermBase[idn[d]]
                        vel = VelWet(xref[j,:],qp,qpp,Ndot,Ndotp)
                        wh[count+c+2*d].dat.data[k][:] = vel[:]
        count = count + 4   
    return wh
\end{lstlisting}

\subsubsection{Discrete Fluid problem}

Corresponding to any velocity field in the rod the induced velocity field in the fluid domain
admits the additive decomposition
\begin{equation}
 \mathbf{u}_h = \mathring{\mathbf{u}}_h + \sum_{i=1}^{N_{\mbox{\tiny dof}}} { \alpha_i \, \mathbf{w}_h^i } = 
 \mathring{\mathbf{u}}_h + \sum_{i=1}^{N_{\mbox{\tiny dof}}} { \alpha_i \,\widehat{H}_{h,\mathbf{q}_H^n}\mathbf{N}^i}
 \label{eq:uhadditive}
\end{equation}
where $\mathring{\mathbf{u}}_h$ vanishes on 
$\mathcal{S}(\mathbf{q}^n_H)$.
Eq. (\ref{eq:uhadditive}) translates into the python code
\begin{lstlisting}[language=Python]
uo,p = split(sol_f) 
u = uo
for i in range(ndofs):
    u += Constant(alpha[i])*wh[i]
\end{lstlisting}
The implementation presented here is restricted 
to the case of $P_1$ linear elements both for velocity and pressure,
rendered stable by taking 
$\tau_h = h_K^2 / 4 \mu$ where $h_K$ stands for the cell diameter. 
In the experiments below no body forces are considered, although their inclusion
is straightforward. The following function determines $\mathring{\mathbf{u}}_h$ and the pressure field $p_h$ by solving the fluid problem (\ref{eq41}):

\begin{lstlisting}[language=Python]
def SolveFluid(u, p, v, q, hk, sol, bcs): # Eq. (45)
    mu = viscosity(u)
    tauh = hk*hk/(4*mu) # Stabilization parameter
    F  = ( inner(2*mu*sym(grad(u)), grad(v)) - div(v)*p )*dx
    F += ( q*div(u) + tauh*inner(grad(q), grad(p)) )*dx
    Jac = derivative(F, sol)
    solve(F == 0, sol, bcs=bcs, J=Jac, solver_parameters=sp)
    uo,p = split(sol)
    return uo,p
\end{lstlisting}
Recall that the fluid problem can be nonlinear if the viscosity $\mu$ 
depends on $\mathbf{u}$. This nonlinearity is
dealt with by means of a Newton-Raphson method with line search.
Otherwise, the fluid solver converges in one single iteration.
Firedrake is coupled to the PETSc library \cite{petsc2021}, from which the snes 
functions for the non-linear part and several solvers for the linear part,
can be invoked. This is defined in the \textcode|solver_parameters| object given next.
Notice, the mumps \cite{mumps2001, mumps2019} direct solver has been chosen for the linear solves.
\begin{lstlisting}[language=Python]
# Fluid-solver options
sp = {
    'snes_type': 'newtonls',
    'snes_monitor': None,
    'snes_converged_reason': None,
    'snes_linesearch_type': 'basic',
    'snes_max_it': 20,
    'snes_divergence_tolerance': 1e20,
    'snes_rtol': 1e-7,
    'snes_atol': 1e-10,
    'mat_type': 'aij',
    'ksp_type': 'preonly',
    'pc_type': 'lu',
    'pc_factor_mat_solver_type': 'mumps'
}
\end{lstlisting}

\subsubsection{FSI nonlinear coupling solver}

Consider Eq. (\ref{eq:residual}) viewed as a function of $\mathbf{\dot{q}}_H^{n+1}$, 
or equivalently as a function of the unknown parameters 
$\underline{\alpha} = (\alpha_1, \dots, \alpha_{N_{\mbox{\tiny dof}}})^{\intercal}$.
For each $i=1,\dots, N_{\mbox{\tiny dof}}$
\begin{equation*}
  \mathcal{R}^i(\underline{\alpha}) =
  \int_{\Omega_f(\mathbf{q}_H^n)} \left [ 2\mu\,\boldsymbol{\varepsilon} (\mathbf{u}_h) \bm{:} \boldsymbol{\varepsilon} (\mathbf{w}_h^i)
-{p}_h\,\nabla \cdot \mathbf{w}_h^i
   \right ] ~d\Omega~
    + %\langle D_{\mathbf{q}}E(t,\mathbf{q}_H^{\star}),\delta \mathbf{q}_H \rangle 
        \,\langle D_{\mathbf{q}}E(t_*,\mathbf{q}_*),\mathbf{N}^i \rangle +
        \end{equation*}
        \begin{equation}
  ~+~ \delta \, t\, S(\mathbf{q}_H^n;\mathbf{\dot{q}}_H^{n+1}%-\mathbf{q}_H^n
  ,\mathbf{N}^i)
  - \,\langle \mathbf{F}(t_*),\mathbf{N}^i \rangle~,
    \label{eq:P2disc}
\end{equation}
%where $\mathbf{F}(t)$ has been assumed to be zero for simplicity.
where now the implementation of all the terms has been made precise.
The residual function \textcode|residualFSI()| has to be provided 
to the nonlinear solver. To find $\underline{\alpha}$ we use 
the function \textcode|root()| from the \textcode|scipy.optimize| 
library (see \cite{2020SciPy-NMeth}) within the loop over time steps as follows
\begin{lstlisting}[language=Python]
def residualFSI(alpha, 
                mesh, hk, v, q, sol, wh, Ndofs, dt, 
                eref, espont, fe, 
                qHn, qHdot, deltaqH, bcs):
    u,p = split(sol)
    for i in range(Ndofs):
        u += Constant(alpha[i])*wh[i]
    
    u,p = SolveFluid(u, p, v, q, hk, sol, bcs)
    for i in range(Ndofs):
        u += Constant(alpha[i])*wh[i]
    
    for i in range(Nnodrod): # Build rod velocity based on alpha's
        qHdot.dat.data[2*i][0]   = alpha[4*i+0]
        qHdot.dat.data[2*i][1]   = alpha[4*i+1]
        qHdot.dat.data[2*i+1][0] = alpha[4*i+2]*h_r
        qHdot.dat.data[2*i+1][1] = alpha[4*i+3]*h_r

    mu = viscosity(u)
    for i in range(Ndofs): # Testing against the wh[:] 
        res_fsi[i] = assemble(inner(2*mu*sym(grad(u)),grad(wh[i]))*dx)
        res_fsi[i]-= assemble(div(wh[i])*p*dx)
        res_fsi[i]+= b_RodTerms_b(dt, eref,espont,fe, qHn,qHdot,deltaqH[i])
    return res_fsi
\end{lstlisting}
\begin{lstlisting}[language=Python]
...
# Begin loop over time steps
for n in range(ntimes):
    qHn.assign(qH)
    espont = SpontaneousStrain(tt, dt, x_r)
    fe = ExternalForce(tt, dt, x_r)
    wh = ApplyHhatqN(W_f, xref, list_nodes_wet, qHn)
    solnl = root(residualFSI,
                 alpha,
                 args    = ((mesh_f, hk, v, q, sol_f, wh, Ndofs, dt, 
                            eref, espont, fe, qHn, qHdot, deltaqH, bcs_f)),
                 method  =  'krylov',
                 options = {'disp': True, 'fatol': 1e-6, 
                            'maxiter': 5, 'line_search': 'armijo',
                            'jac_options': {'inner_tol': 1e-12, 
                                            'inner_maxiter': 100000, 
                                            'method': 'gmres', 
                                            'rdiff': 1e-4}})

    alpha = solnl.x
    ...
    qH.assign(qHn + Constant(dt)*qHdot)
    new_wet = UpdateWetBoundary(xref, qH)
    for i in range(nnod_wet): # Move the wet boundary
        x_f[list_nodes_wet[i],:] = new_wet[i,:]
    
    if(need_remesh):
        new_mesh_f, ... = RemeshFluidRegion ( mesh_f, ...)
...

\end{lstlisting}
In the snippet above, the code has been slightly simplified so as to 
point out to the important features of the implementation. 
The solution is sought by using a Newton-Krylov method \cite{newton_krylov}.
Being an inexact Newton method, convergence speed may depend 
on the choice of algorithmic parameters. 
In the numerical experiments below, the set of parameters shown in the code have worked well in general and can be recommended. Finally, results can be visualized with 
Paraview \cite{Paraview} as the simulation runs simply by opening some \textcode|pvd| 
files created to output the velocity and pressure fields. In the case of a nonlinear 
fluid an additional file to output the viscosity is also saved.
Additional information, such as the rod elastic energy, the fluid dissipation
and the geometric position of the left end of the rod as a function of time
are written to an ascii file called \textcode|evol.txt|.

The source code that solves the last numerical experiment to be presented
in the next section, namely, the swimming of a flexible rod immersed in a viscous fluid
is provided. The code can be found in \textcode|https://gitlab.com/rfausas/microswimmers.git| 
along with details on the exact firedrake version, its dependencies used in this work 
and some useful information to run the script (see \textcode|README|).  
All experiments to be presented hereafter were run into a desktop computer 
with Core i7-7700K CPU @ 4.2GHz with Linux Ubuntu 18.04.3 LTS.

\section{Numerical experiments}

\subsection{Roll-up and roll-out  of a flexible rod}

\begin{figure} 
   \begin{center}
       \scalebox{0.9}{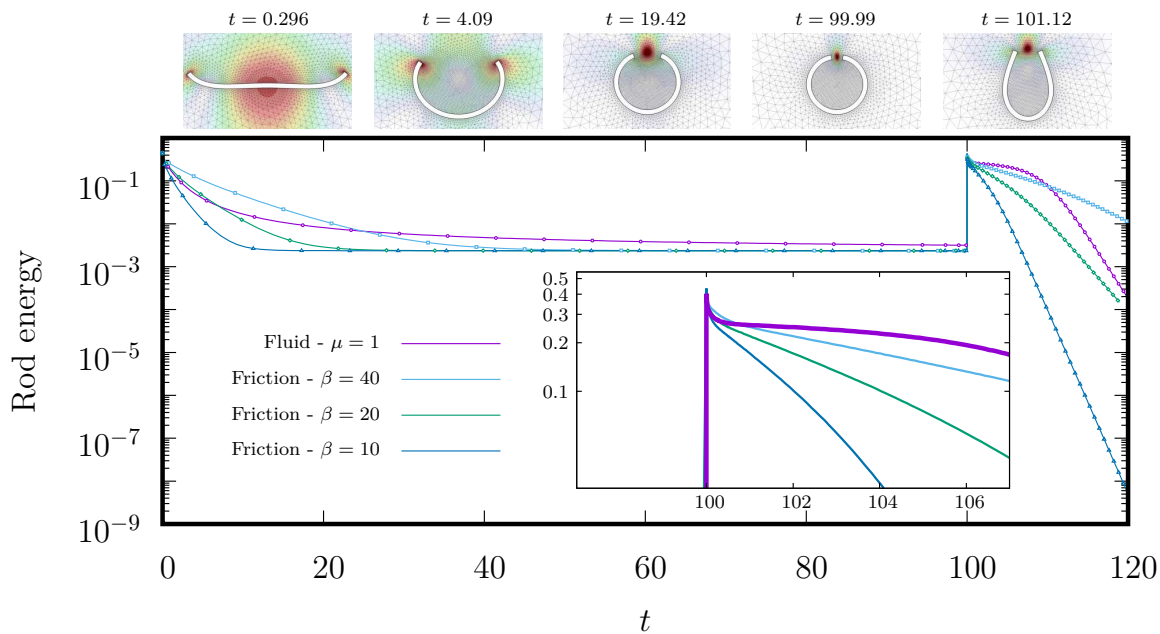}
       \caption{Energy as a function of time for the roll-up and release of a pre-stressed rod. The purple curve corresponds to the fluid-structure-interaction model proposed in this article, with the fluid being Newtonian with $\mu = 1$. The other curves consider a simpler model, in which the fluid is replaced by a friction force $-\beta \mathbf{u}$ with $\beta = 10,~20$ and $40$. The inserts show the geometry and fluid mesh at several instants, colored with the velocity magnitude.
       }
       \label{fig:roll_up_evol_closure}
   \end{center}
\end{figure}

\begin{figure} 
   \begin{center}
       \scalebox{0.095}{\input{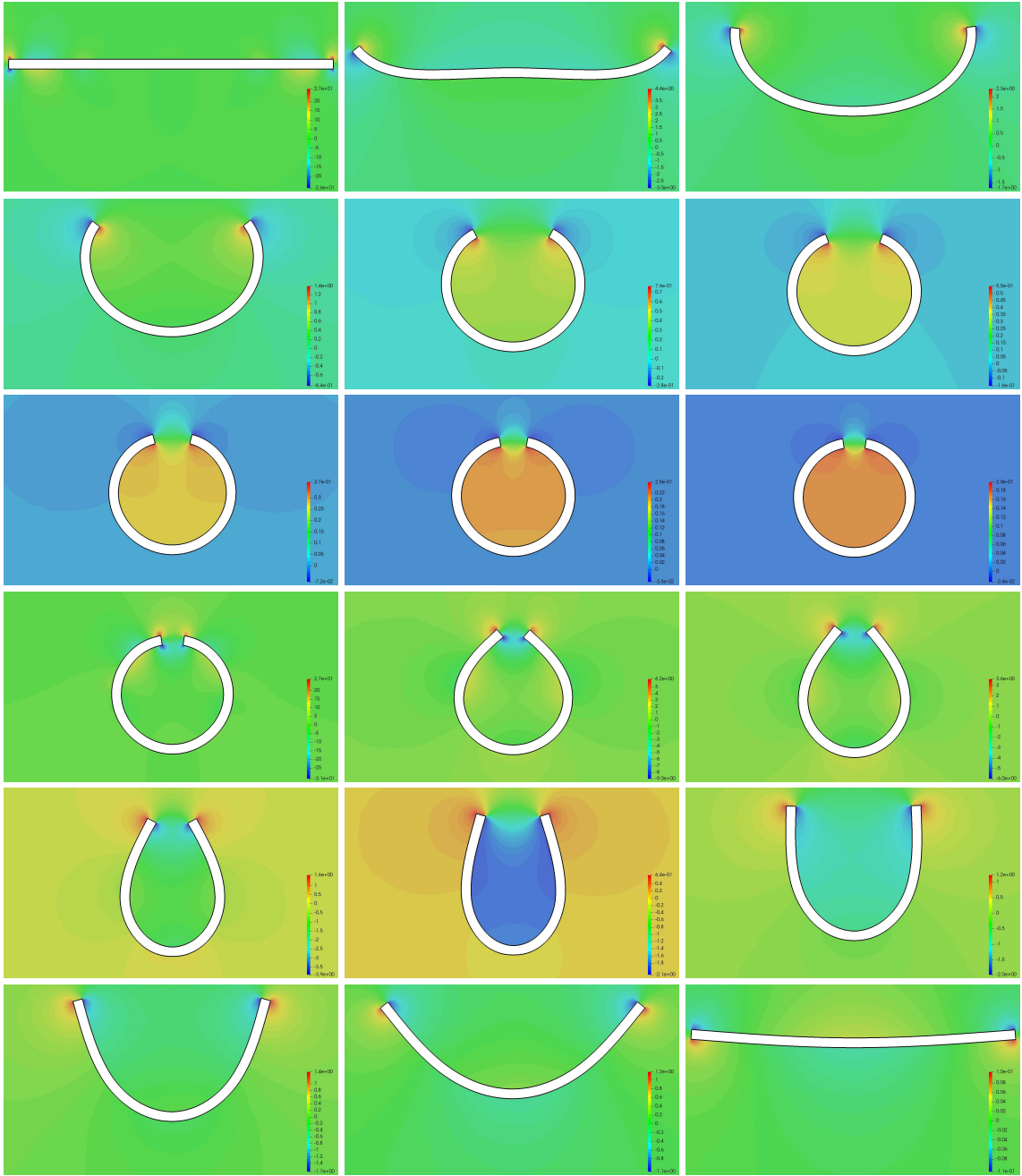}}
       \caption{Shape evolution and pressure field for the roll-up and release of a pre-stressed rod immersed in a viscous fluid. The instants of the frames above are 
       $t=0,~0.296,~1.698,~4.089,~11.42,~19.42,~37.42,~67.42,~99.99,~100,100.07,~100.2,~100.5,~103.0,$  $~107.2,~109.0,~111.0,~119.0$.} %$t=0,~0.296,~1.698,~4.089,~11.422,~19.425,~37.425,~67.425,99.99,~100,100.07,~100.2,~100.520,~103.027,~107.27,~109.027,~111.027,~119.027$.}
       \label{fig:frames}
   \end{center}
\end{figure}

The first experiment we consider consists of the roll-up and the subsequent roll-out of an initially straight rod of unit length.
This is accomplished by setting the spontaneous curvature 
$\kappa_0$ to
\begin{equation}
\kappa_0(s,t) = \left \{
\begin{array}{lr}
-0.999 \,\times\, 2\pi & \mbox{if}~t < 100\\
0 & \mbox{otherwise}
\end{array}
\right.
\end{equation}
i.e., the initially straight rod rolls towards a circular rod 
of radius $R = \frac{1}{0.999 \, 2\pi }$ (i.e., almost to
closure) and at $t = 100$ the rod is released. 
We take $C_{\epsilon} = 90$ and $C_{\kappa} = 0.0225$ and
$a(s) = \frac{1}{\parallel \textbf{q}_H'(s) \parallel}$.
Both ends of the rod are flat. 
The computational domain is the square region $\Omega = [0,3]^2$.
Both lateral walls are open (traction free) and the bottom and top walls are subject to zero-velocity
boundary conditions. The rod is discretized with $N_{\mbox{\tiny el}} = 4$ Hermite elements
and the fluid mesh has elements of characteristic size $h = 1/100$
close to the rod's boundary. Both stages of the process are
shown in Fig. \ref{fig:roll_up_evol_closure} that displays
the rod energy as a function of time. Initially, the rod
relaxes to a minimum energy state ($\mathcal{E} \sim 2.3 \times 10^{-3}$). 
This minimum is not zero because the circular shape can not be represented
exactly by cubic polynomials. At $t = 100$ the pre-stress
is released and a second relaxation process takes place until the initial straight shape is recovered, exactly. 
For this problem the time step is adaptively chosen
so as to better capture the different time scales
observed along the simulation, specially in the initial stages ($t < 1$)
and immediately after release $t \ge 100$.
Inserts with plots of the rod shape and the triangular fluid mesh generated are shown at different times, including the rod
just prior to release and just after it. In Fig. \ref{fig:roll_up_evol_closure}
we have also included the behavior of a simpler model of rod, in which the ambient fluid
is substituted by a friction force $-\beta \, \mathbf{u}_h$
proportional to the rod's surface velocity. To be precise, the variational residual of this simpler model reads 
%(assuming zero external forcing)
\begin{equation}
  \mathcal{R}^i(\underline{\alpha}) =
  \int_{\mathcal{S}(\mathbf{q}_H^n)} \beta \, \mathbf{u}_h \cdot \mathbf{w}_h^i ~d\mathcal{S}~
    + \langle D_{\mathbf{q}}E(t_*,\mathbf{q}_*),\mathbf{N}^i \rangle + \delta t S(\mathbf{q}_H^n,\mathbf{\dot{q}}_H^{n+1},\mathbf{N}^i) - \langle \mathbf{F}(t_*),\mathbf{N}^i \rangle = 0~.
\end{equation}
The plot shows results corresponding to $\beta = 10,~20$ and $40$.
One thing to notice is that the dynamical behavior of the simpler model differs
from that of the true fluid structure interaction problem,
specially in the released phase. This shows the importance
of the complete FSI treatment being proposed in this work. 
Snapshots of the complete rod evolution together with contours of pressure
are shown in Fig. \ref{fig:frames}, the first three rows
of frames corresponding to the pre-stressed phase and the last 
three to the released phase. The dynamics and intermediate
shapes of these two phases are appreciably different. 

\begin{figure} 
   \begin{center}
       \scalebox{0.65}{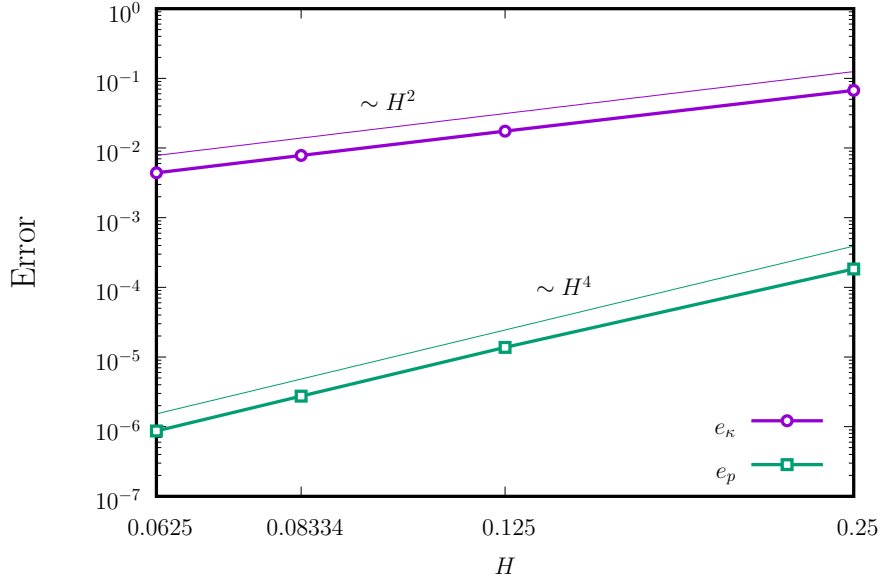}
       \caption{Curvature and positional error for the roll-up of a pre-stressed rod
       immerse in a viscous fluid.}
       \label{fig:error_vs_H}
   \end{center}
\end{figure}

\begin{figure} 
   \begin{center}
       \scalebox{0.65}{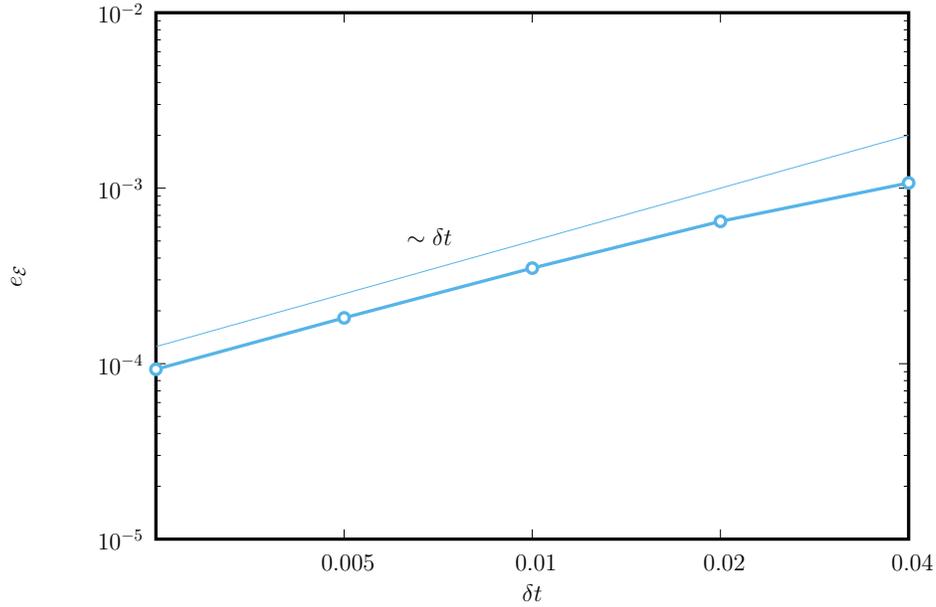}
       \caption{Energy balance error for the roll-up of a pre-stressed rod immersed in a viscous fluid.}
       \label{fig:error_vs_dt}
   \end{center}
\end{figure}

Let us now assess the convergence properties of the scheme by pre-stressing the rod with $\kappa_0 = -\pi$, such that equilibrium shape corresponds to half of a circle of radius $1/\pi$. We compute the curvature and the positional errors as
\begin{equation}
e_{\kappa} =  \left ( \int_{0}^1{\kappa^2}\,ds \right )^2 , ~~~~e_p = \int_{0}^1{\left |  (\tilde{\mathbf{q}}_H\cdot \tilde{\mathbf{q}}_H)^{\frac12}  - \pi^{-1} \right | }\,ds
\end{equation}
with $\tilde{\mathbf{q}}_H = {\mathbf{q}}_H - \mathbf{r}_c$ 
being $\mathbf{r}_c = \frac12 \left ( \mathbf{q}_H(s=0) + \mathbf{q}_H(s=1) \right )$.
Fig. \ref{fig:error_vs_H} shows these errors as a function of the 
mesh parameter $H$ at the final time of the relaxation.
The convergence rate for the curvature error $e_{\kappa}$ is $\mathcal{O}(H^2)$ 
and for the positional error $e_p$ is $\mathcal{O}(H^4)$, as expected for Hermite cubic polynomials \cite{Gebhardt2021}. 
To conclude this subsection, the convergence properties of the semi-implicit scheme
with respect to time step refinement is also assessed. To that end, we compute
a discrete energy balance error as 
\begin{equation}
 e_{\mathcal{E}} = \left | \int_{t_i}^{t_f}{R_{\mathcal{E}}(t)}dt\right | \label{eq:L2t}
\end{equation}
in which
\begin{equation}
 R_{\mathcal{E}}(t) = \frac{d\mathcal{E}}{dt} +
 \int_{\Omega(\mathbf{q}_H)}{\left ( 2\mu\nabla^S{\bf u}_h : \nabla^S{\bf u}_h - p_h \nabla\cdot{\bf u}_h \right ) \,d{\Omega}} ~+~ \delta{t} \, S(\mathbf{q}_H;\mathbf{\dot{q}}_H,\mathbf{\dot{q}}_H) \label{eq:deltaE_balance}
\end{equation}
is the discrete power imbalance, namely the rate of energy change (first term of $R_{\mathcal{E}}$) plus fluid dissipation (second term) plus the dissipation produced by the stabilization (third term).  
Notice that the $p_h \nabla\cdot{\bf u}_h$ term in the instantaneous energy balance (\ref{eq:deltaE_balance}) is zero in the exact problem and also in the discrete problem if $\tau_h=0$, but not in our $P_1/P_1$ implementation. In computing (\ref{eq:deltaE_balance}), $\frac{d\mathcal{E}}{dt}$ is taken as $\frac{\mathcal{E}(t_{n}) - \mathcal{E}(t_{n-1})}{\delta{t}}$ and all other quantities are taken at instant $n$. Also, $t_i = 0.04$
and $t_f = 0.8$. The error $e_{\mathcal{E}}$ for different
time steps $\delta{t}$, namely, $0.04,~0.02,~0.01,~0.005$ and $0.0025$ is displayed in Fig. \ref{fig:error_vs_dt},
showing the expected convergence rate $\mathcal{O}(\delta{t})$.
To quantify the relative importance of the different contributions to the energy 
balance (\ref{eq:deltaE_balance}) we show in Table \ref{tab:termsdE} the integral of each
term in the same time interval, $(t_i,t_f)$. The main contributions are the physical ones, namely the viscous dissipation and the time derivative of $\mathcal{E}$ (second
and fifth columns) which approach each other as the time step is refined. The semi-implicit stabilization term (fourth column) is smaller and $\mathcal{O}(\delta t)$. The pressure term (third column) is also rather small but independent of $\delta t$, since in fact it is a consequence of the discretization of the fluid problem and thus depends on $h$.
%Notice that $\int_{t_i}^{t_f}{\frac{d\mathcal{E}}{dt}}\,dt$ is simply given by
%$\mathcal{E}(t_f) - \mathcal{E}(t_i)$.
\begin{table}[ht!!]
\begin{center}
\caption{Contributions of the different terms in the discrete energy balance.} \label{tab:termsdE}
~\\
\begin{tabular}{c|c|c|c|c}
$\delta{t}$ &  $\int_{t_i}^{t_f}  \int_{\Omega}{2\mu\nabla^S{\bf u}_h : \nabla^S{\bf u}_h }\,d{\Omega} \,dt$  &  $-\int_{t_i}^{t_f} \int_{\Omega}{p_h\nabla\cdot{\bf u}_h }\,d{\Omega} \,dt$ &  $\int_{t_i}^{t_f} {\delta{t}\,S(\cdot;\cdot,\cdot)}\,dt$ & $\mathcal{E}(t_f) - \mathcal{E}(t_i)$ \\
\hline
0.04 &  $0.037543$  &  $0.00103$ &  $0.003307$ & $-0.042951$ \\
0.02 &  $0.038663$  &  $0.00105$ &  $0.001369$ & $-0.041734$ \\
0.01 &  $0.039334$  &  $0.00107$ &  $0.000702$ & $-0.041462$ \\
0.005 & $0.039705$  &  $0.00108$ &  $0.000359$ & $-0.041335$ \\
0.0025 & $0.039905$  &  $0.00109$ &  $0.000182$ & $-0.041273$
\end{tabular}
\end{center}
\end{table}

\subsection{Large deflection of a flexible beam in a channel}

In the second experiment we consider a cantilever beam of initial
length $S=1$ and thickness $e = 0.03$ 
whose left end is clamped at $(x_1,x_2) = \left (0,\frac{3}{2}\right )$ 
in the computational domain $[0,3]^2$. 
The boundary conditions at the left end of the rod are thus $\mathbf{q}=\left ( 0,\frac32\right )$ and $\frac{dq_2}{ds}=0$. 
The right end of the beam is rounded and free. 
The material parameters are $C_\epsilon = 90$ 
and $C_\kappa = 0.0225$ as in the previous experiments. 
The fluid is assumed to be Newtonian with viscosity
$\mu = 1$. Also, we consider the function $a(s) = 1$ (see Eq. (\ref{eq:aofs})).
The beam is discretized with $N_{\mbox{\tiny el}} = 8$ elements and the fluid mesh 
with elements having characteristic size $h = 1/100$ near the beam boundary,
which gives meshes with varying number of elements, that oscillates between $\sim5$K to 
$\sim6.5$K. The time step $\delta{t}$ is set to the constant value $0.05$ throughout the simulation.
We impose a non-uniform (constant in time) spontaneous curvature given by
\begin{equation}
 \kappa_0(s) =  2\pi \left( 1 + \tanh(4s) \right).
\end{equation}
As a result of this, the beam starts to deform into a scroll shape from $t = 0$ until
$t = 12$ when a fluid inlet located on the lower part of the 
left boundary is suddenly opened. The fluid velocity profile at the inlet is set to be
the function
\begin{equation}
\mathbf{u}(x_1=0,x_2) = \left ( (3 - e)\,x_2^3 - 2 x_2^4,~0 \right),~~~~0 \le x_2 \le \frac{3-e}{2}~.
\end{equation}
The influx of the fluid, which leaves through the upper boundary, unrolls the beam and pushes it upwards, resulting in a complex process with large deformations and rotations and strong fluid-beam interaction.
The bottom wall ($x_2 = 0$), the right wall ($x_1 = 3$) and 
the part of the left wall above the beam ($x_1=0, x_2 > \frac{3+e}{2}$), are permanently subject to zero-velocity boundary conditions. The top wall ($x_2 = 3$) is an outlet boundary through which the fluid is free to leave 
the domain.

\begin{figure} 
   \begin{center}
       \scalebox{0.093}{\input{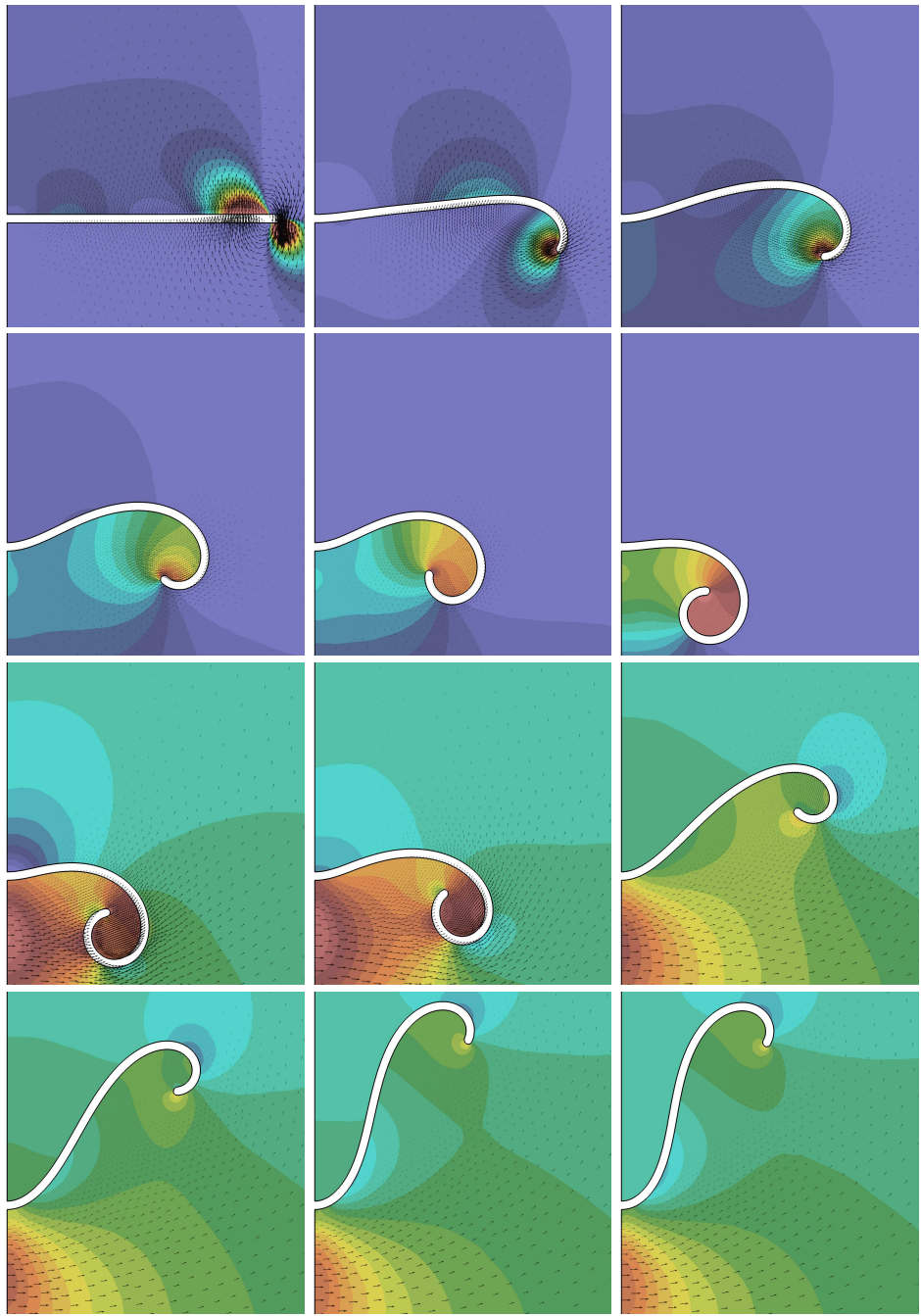}}
       \caption{Shape, contours of pressure field and velocity vectors 
         for the cantilever rod subject to a pre-stress and and a velocity field switched on from time $t = 12$.
       Pressure scale: $4.1$ (red), $-0.11$ (blue) for $t < 12$ (first two rows) 
       and $9.7$ (red), $-4.7$ (blue) for $t \ge 12$ (last two rows). 
       Velocity vector scale: $0.15$ ($t < 12$) and $0.035$ ($t \ge 12$). 
       Shown are the instants 
       $t=0,~0.5,~1.25,~2.5,~5,~11.95,~12.1,~12.4,~14,~16,~25,~30$.}
       \label{fig:cantilever_pres}
   \end{center}
\end{figure}

\begin{figure} 
   \begin{center}
      \scalebox{0.1475}{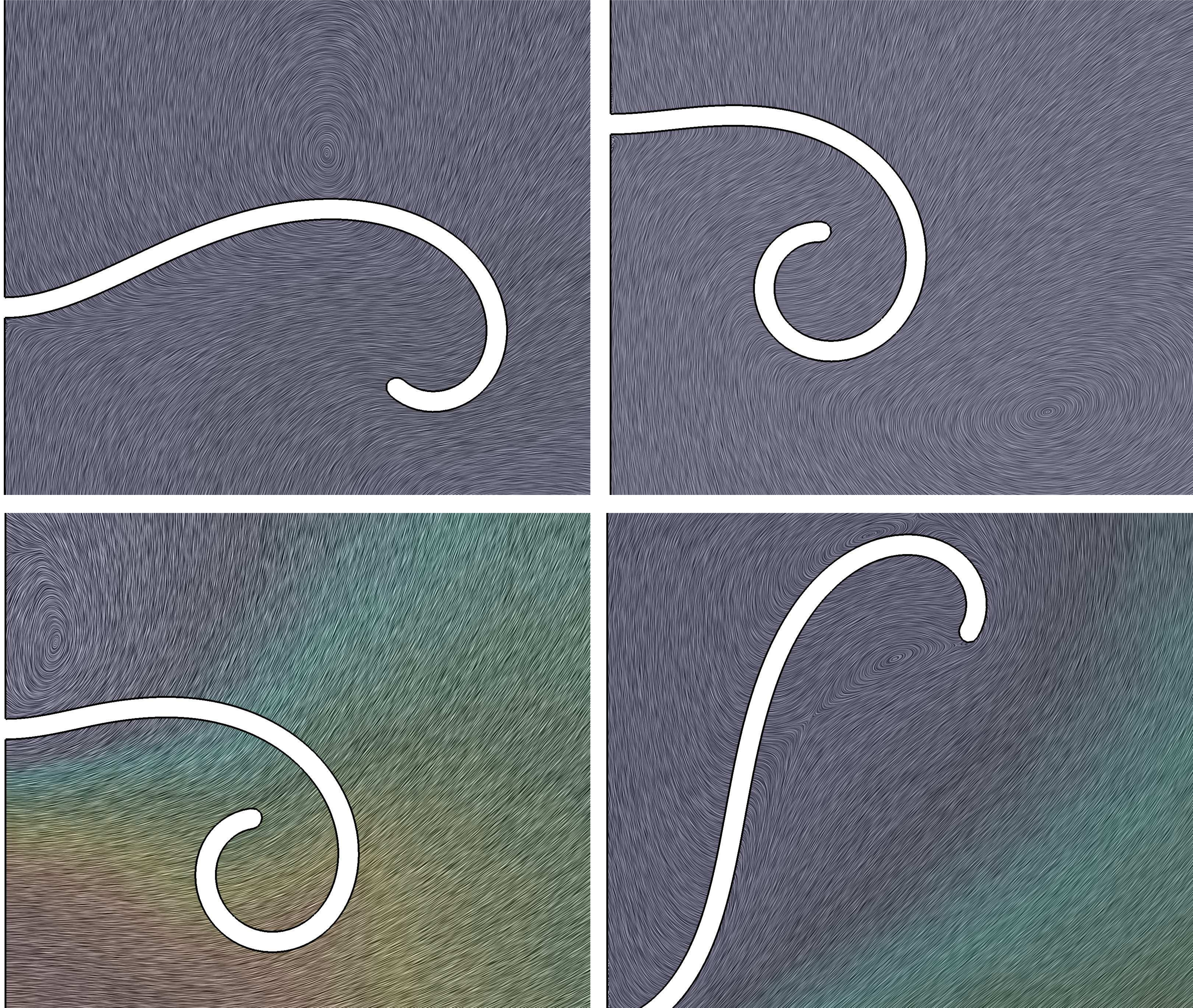}
       %\scalebox{0.093}{\input{figures/new_cantilever_frames.pdf_t}}
       \caption{LIC Streamlines colored by the velocity magnitude
       for the cantilever rod subject
       to a pre-stress and a velocity field switched on at time $t = 12$. 
       Velocity scale: $1.2$ (red), $0.3$ (cyan) and $0$ (dark blue).
       Shown are details corresponding to instants $t=2.5,~11.95,~12.1,~30$.}
       \label{fig:cantilever_vel}
   \end{center}
\end{figure}

The complete deformation process is plotted in Fig. \ref{fig:cantilever_pres}
that shows the pressure field and arrows of the fluid velocity.
Due to the large difference in magnitude, the velocity arrows in the roll-up
phase (first two rows of frames) are scaled by the factor $0.15$,
while in the roll-out phase (the last two rows of frames) this factor is $0.035$.   
The beam attains a stationary equilibrium shape at the final time of the simulation
as the result of the balance between the internal beam forces and the interaction with the ambient fluid.
Also, Fig. \ref{fig:cantilever_vel} displays details with LIC (Line Convolution Integral) 
streamlines produced by Paraview 5.4.1, colored by the velocity magnitude 
at some selected instants.

\subsection{Swimming of a finite-length flexible rod}

In this last experiment a swimming flexible rod immersed in a viscous fluid
is studied. The computational domain is the square region $\Omega = [0,3]^2$.
Both lateral walls are open and the bottom and top wall are subject to zero-velocity
boundary conditions.
The rod shown in Fig. \ref{fig:refswimmer} is considered, with $e = 0.03$, $D_c = 3\,e$, $H_c = 3\,e$, and unit initial length. The elastic parameters are taken to be $C_\epsilon = 90$ 
and $C_\kappa = 0.0225$. 

\begin{figure} 
   \begin{center}
       \scalebox{0.15}{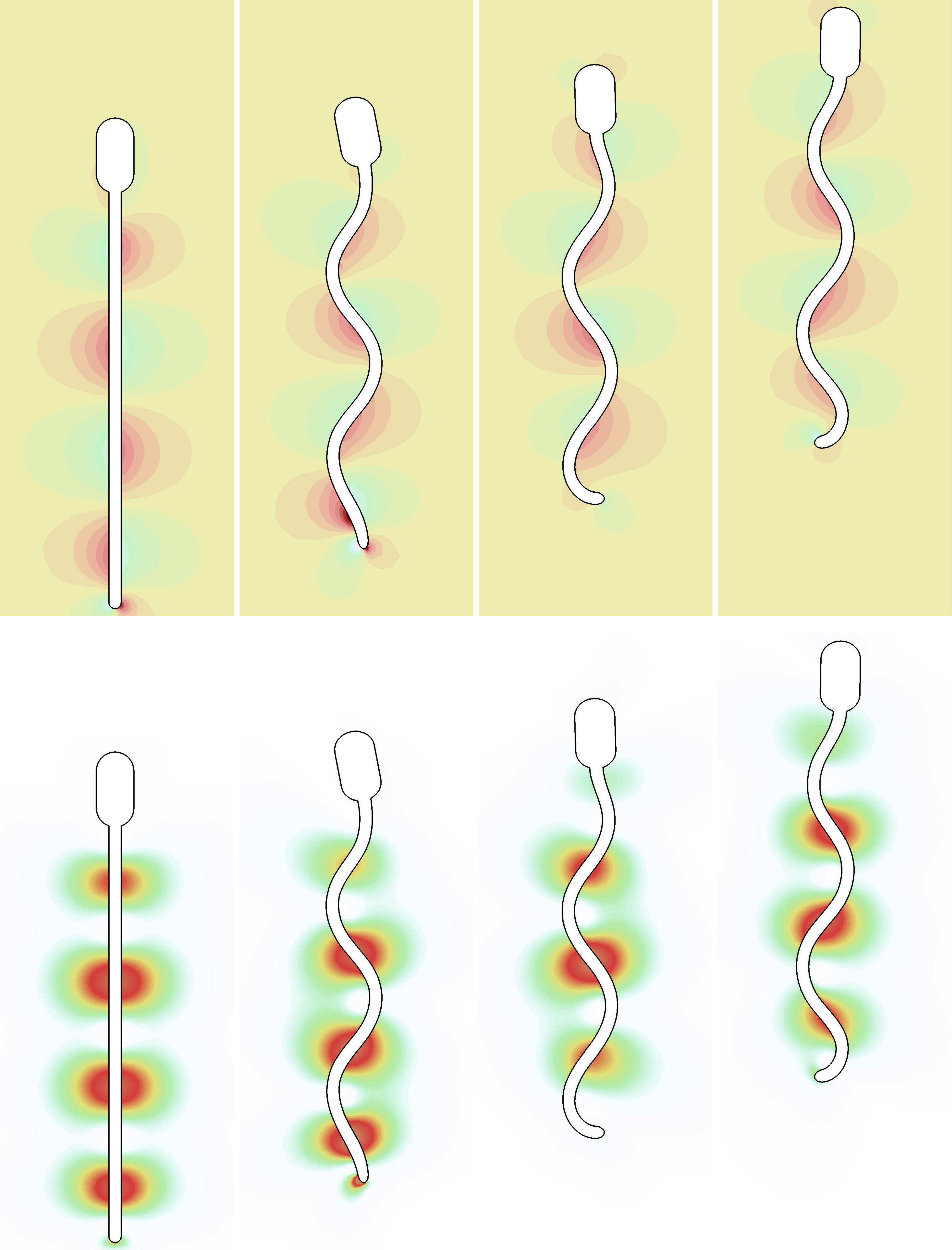}
       \caption{Contours of velocity magnitude (top) and pressure (bottom) for a swimmer in a Newtonian fluid at different times (from left to right, 
       $t = 0.0,~1.0,~2.4,~3.2$). The rod is discretized with $H = 1/8$ and 
       the fluid with typical refinement $h = 1/100$.}\label{fig:swimmer_frames_base}
   \end{center}
\end{figure}

\begin{figure} 
   \begin{center}
       \scalebox{0.1}{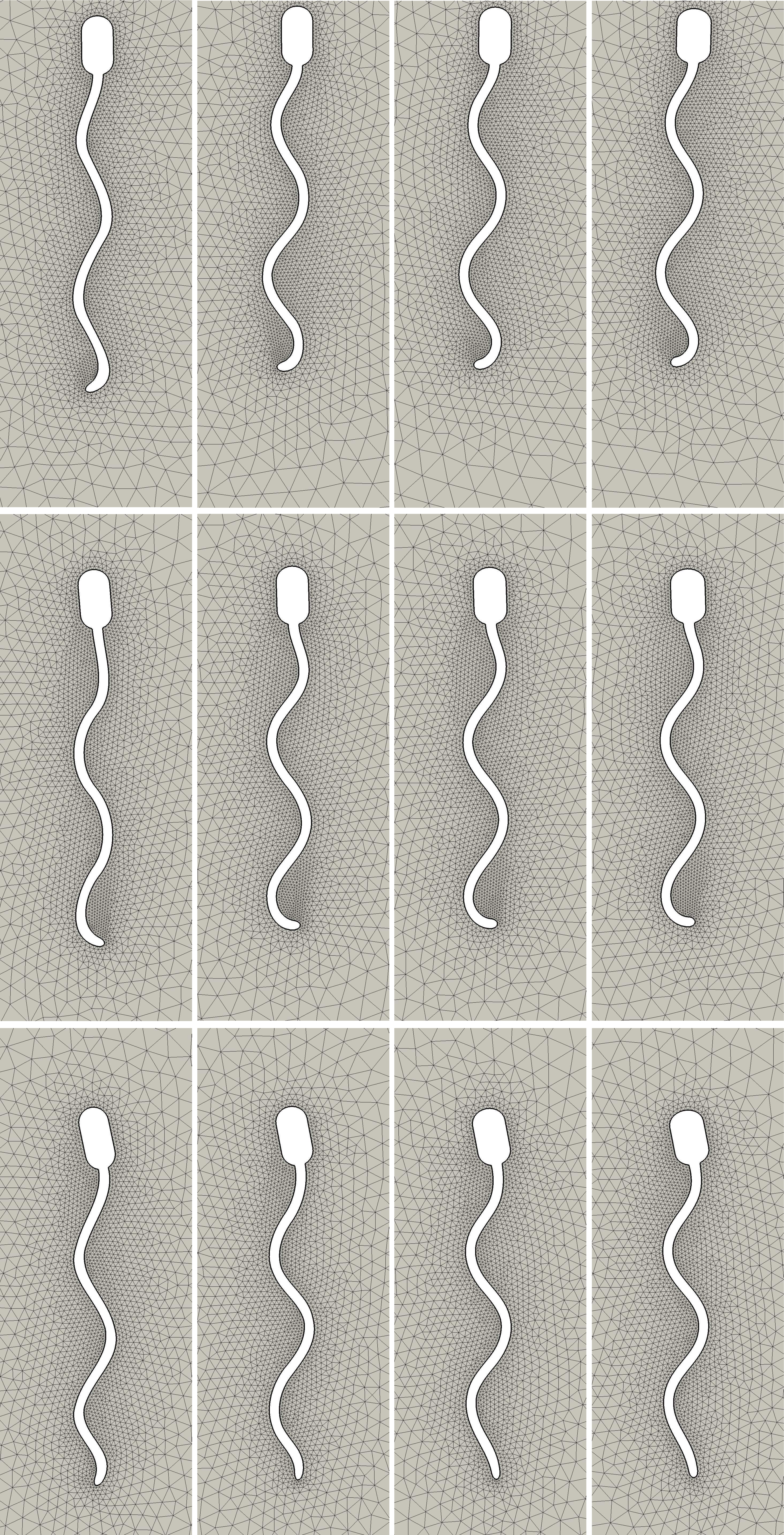}
       \caption{Fluid mesh for a swimmer in a Newtonian fluid
       with typical refinement $h = 1/100$ at different
       times (from bottom to top, $t = 1.0,~1.8,~3.2$). The rod is
       discretized with different number of elements (from left to right, $H = 1/4,~1/6,~1/8,~1/12$).}
       \label{fig:swimmer_effec_rodmesh}
   \end{center}
\end{figure}

\begin{figure} 
   \begin{center}
       \scalebox{0.825}{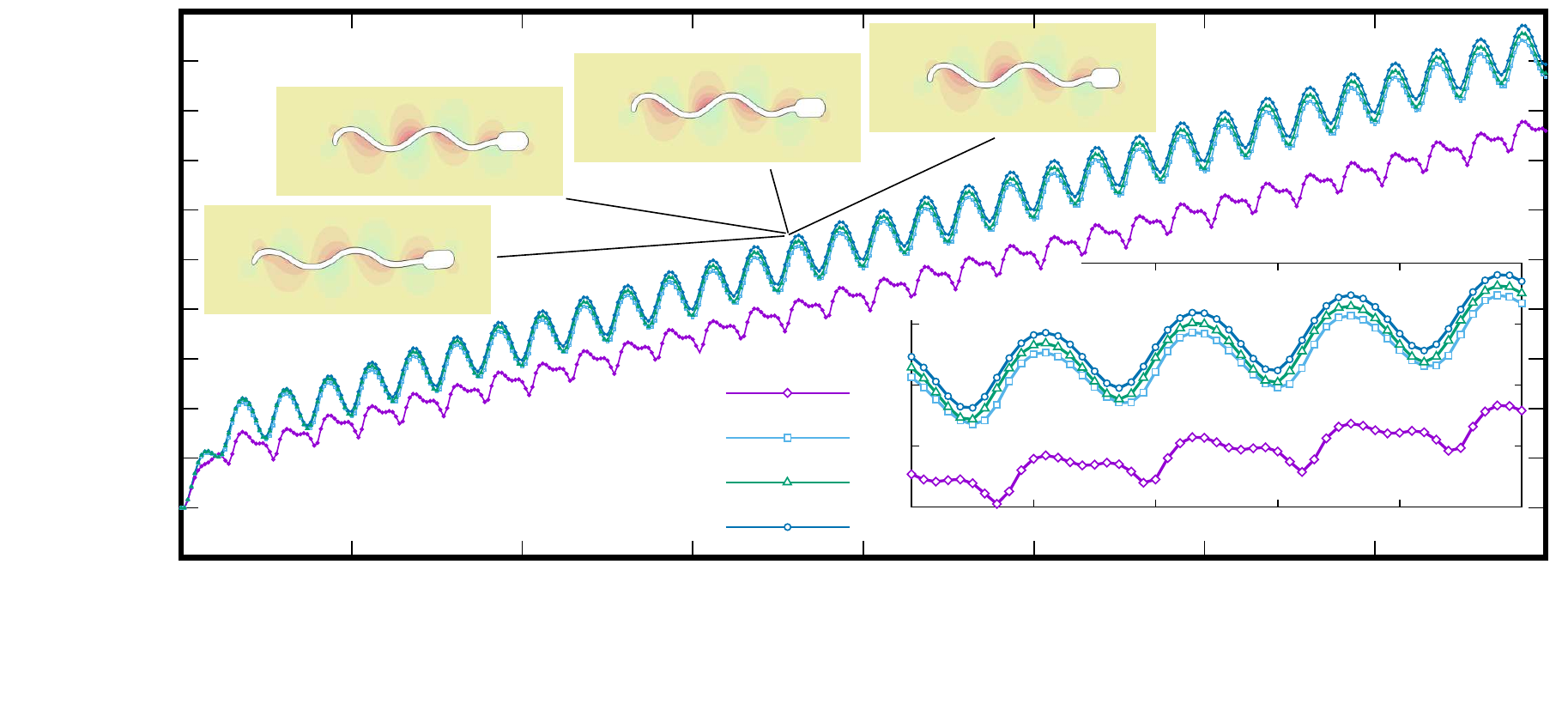}
       \caption{Swimmer position as a function of time.
       The fluid is discretized with mesh refinement 
       $1/100$. The rod is discretized with different number
       of elements $H = 1/4$ (purple), $H = 1/6$ (light blue),
       $H = 1/8$ (green), $H = 1/12$ (dark blue).}\label{fig:evol_swimmer_effect_rodmesh}
   \end{center}
\end{figure}

\begin{figure} 
   \begin{center}
       \scalebox{0.1}{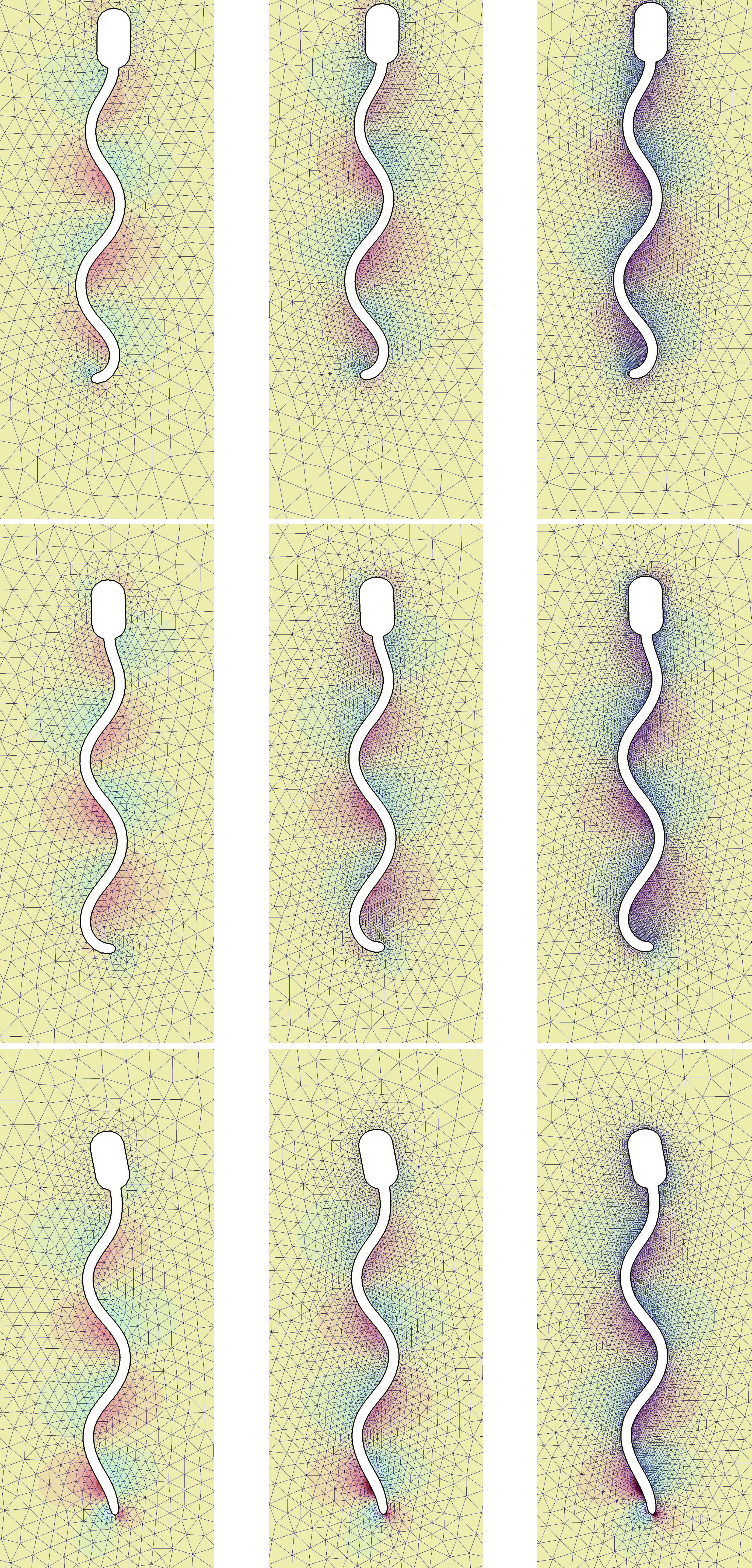}
       \caption{Contours of pressure and mesh for a swimmer 
       in a Newtonian fluid at different times (from bottom to top,  
       $t = 1.0,~1.8,~3.2$). The rod is discretized with $H = 1/8$ and 
       the fluid with different mesh refinements 
       (from left to right, $h = 1/50,~1/100,~1/200$).}\label{fig:swimmer_effect_fluidmesh}
   \end{center}
\end{figure}

\begin{figure} 
   \begin{center}
       \scalebox{0.825}{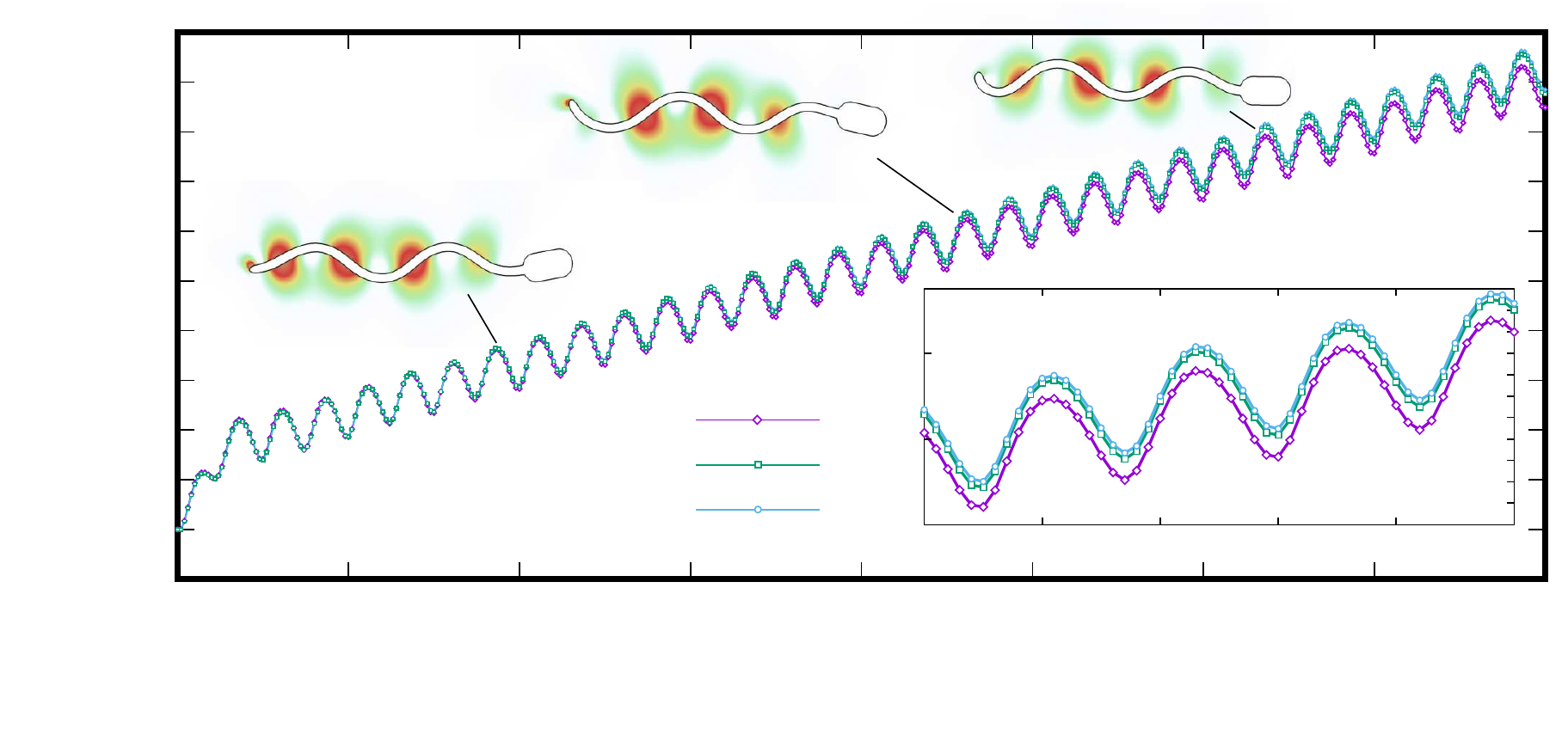}
       \caption{Swimmer position as a function of time. The rod 
       is discretized with $H = 1/8$ and the fluid
       is discretized with diffrent mesh refinements
       $h = 1/50$ (purple), $1/100$ (green) and $1/200$ (blue).}\label{fig:evol_swimmer_effect_fluidmesh}
   \end{center}
\end{figure}

First, the fluid is assumed to be Newtonian with viscosity $\mu = 1$.
By taking the spontaneous curvature in Eq. (\ref{eq:kappa}) to be the wave
\begin{equation}
 \kappa_0(s,t) = 20 \, \sin\left(4\,\pi(s - 2\,t)\right)~,
\end{equation}
the rod swims as seen in Fig. 
\ref{fig:swimmer_frames_base}, which shows the rod's shape together with contours of
pressure and velocity magnitude at different times,
namely, $t = 0.0,~1.0,~2.4,~3.2$. Note that by setting
$\kappa_0(s,t)$ we are not imposing the kinematics of the
rod. Instead, the rod's configurations result from the interaction with the ambient fluid
as the rod tends to relax to its instantaneous minimum energy state.
The time step adopted for all cases in this experiment is $\delta{t} = 0.01$, 
constant throughout the simulation.
The rod is discretized by $8$ elements and the fluid by a triangular mesh with
characteristic element size $h = 0.01$ close to 
the swimmer's boundary. Notice that far away from the swimmer the mesh is coarser.
Remeshing of the fluid domain is being executed at each time step, the total number of 
elements in the triangulation varying around $5$K. This is almost mandatory due to the 
large time step being used in the simulation. We observe that as a result of this periodic
deformation with large amplitude the rod exhibits a net displacement.
As previously stated we consider $a(s) = \frac{1}{\parallel \textbf{q}_H'(s) \parallel}$
(see Eq. (\ref{eq:aofs})), which penalizes to some degree extensional deformations
of the swimmer. Recall that, according to Eq. (\ref{eq:sigmapart}) all points belonging to 
the swimmer's head are attached to the end point of the rod making the head to be
aligned with the rod axis at all times.
To assess the impact of the discretization parameters on the rod's swimming performance,
Fig. \ref{fig:swimmer_effec_rodmesh} shows the swimmer
and the fluid mesh at times $t = 1.0,~1.8,~2.4,~3.2$ (from bottom to top),
and rod resolutions $H = 1/4,~1/6,~1/8,~1/12$ (from left to right).
It is interesting to note that the swimmer's head experiences
some stretching, specially for the coarsest level $H=1/4$.
The level of discretization also affects the net displacement of the swimmer as shown in Fig. \ref{fig:evol_swimmer_effect_rodmesh}.
The main conclusion that can be drawn from these results is that 
the deformation of the rod is surely well represented by a mesh consisting of at least $8$ elements, the results being quite similar for finer discretizations.
Next, Fig. \ref{fig:swimmer_effect_fluidmesh} shows results for different fluid mesh refinements close to the wet boundary of the deformable rod, namely, from left to right, $h = 1/50,~1/100,~1/200$ (corresponding to fluid meshes with about $3$K, $5$K and $10$K elements, respectively). For the rod we consider a fixed resolution $H = 1/8$. Contours of pressure and the fluid 
mesh are shown. The net displacement of the rod 
for this case is plotted in Fig. \ref{fig:evol_swimmer_effect_fluidmesh}. As one can notice
the results are almost independent of the fluid discretization for meshes
finer than $h=1/100$.

Let us take advantage of this test case to compare the different options for time discretization that have been 
introduced in the previous section, namely, the semi-implicit, the explicit and the pseudo-implicit schemes. In all previous tests we have only shown the results of the semi-implicit method, which is in fact the one with the best performance, but it is instructive to report at least one quantitative comparison.  

\begin{figure} 
   \begin{center}
       \scalebox{0.8}{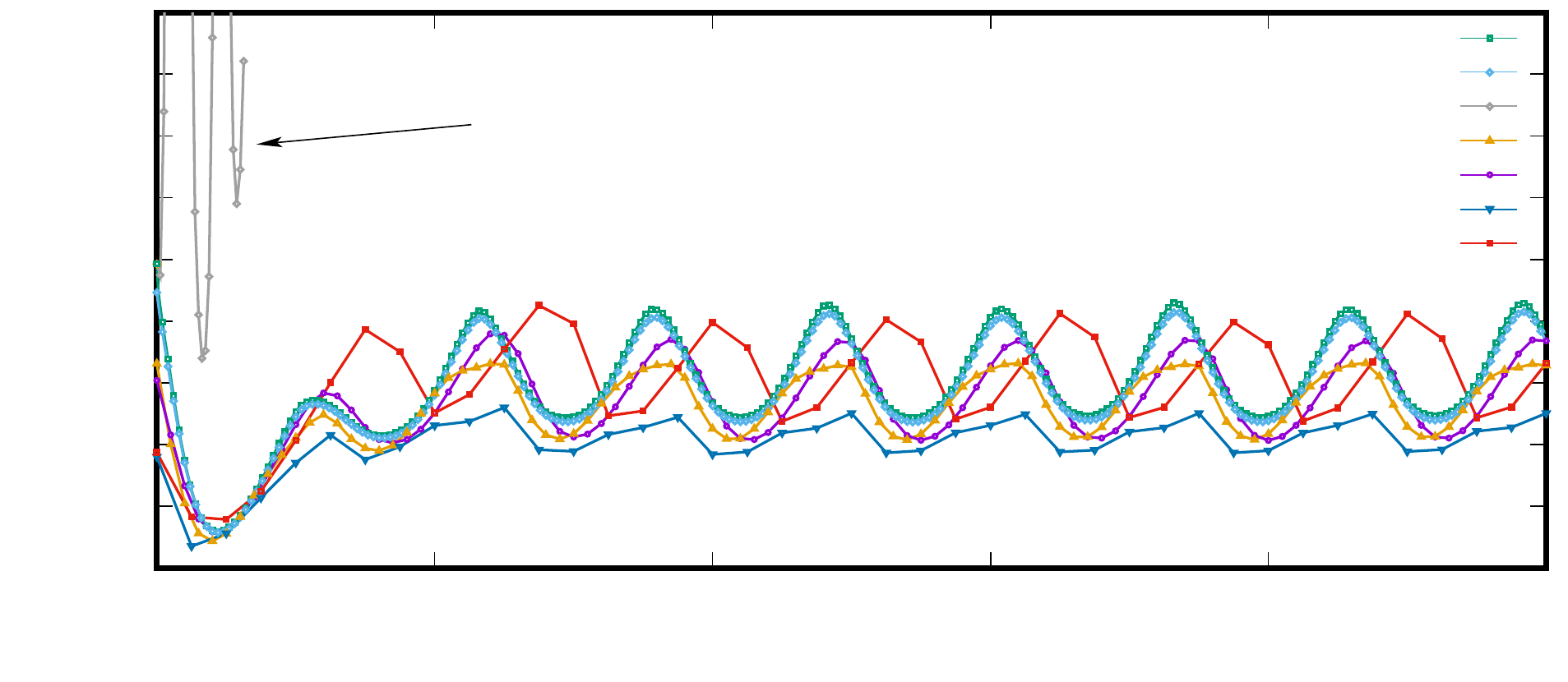}
       \caption{Fluid dissipation as a function of time for the semi-implicit,
       explicit and pseudo-implicit schemes and different time steps.}\label{fig:compa_dissip_methods}
   \end{center}
\end{figure}

\begin{figure} 
   \begin{center}
       \scalebox{0.8}{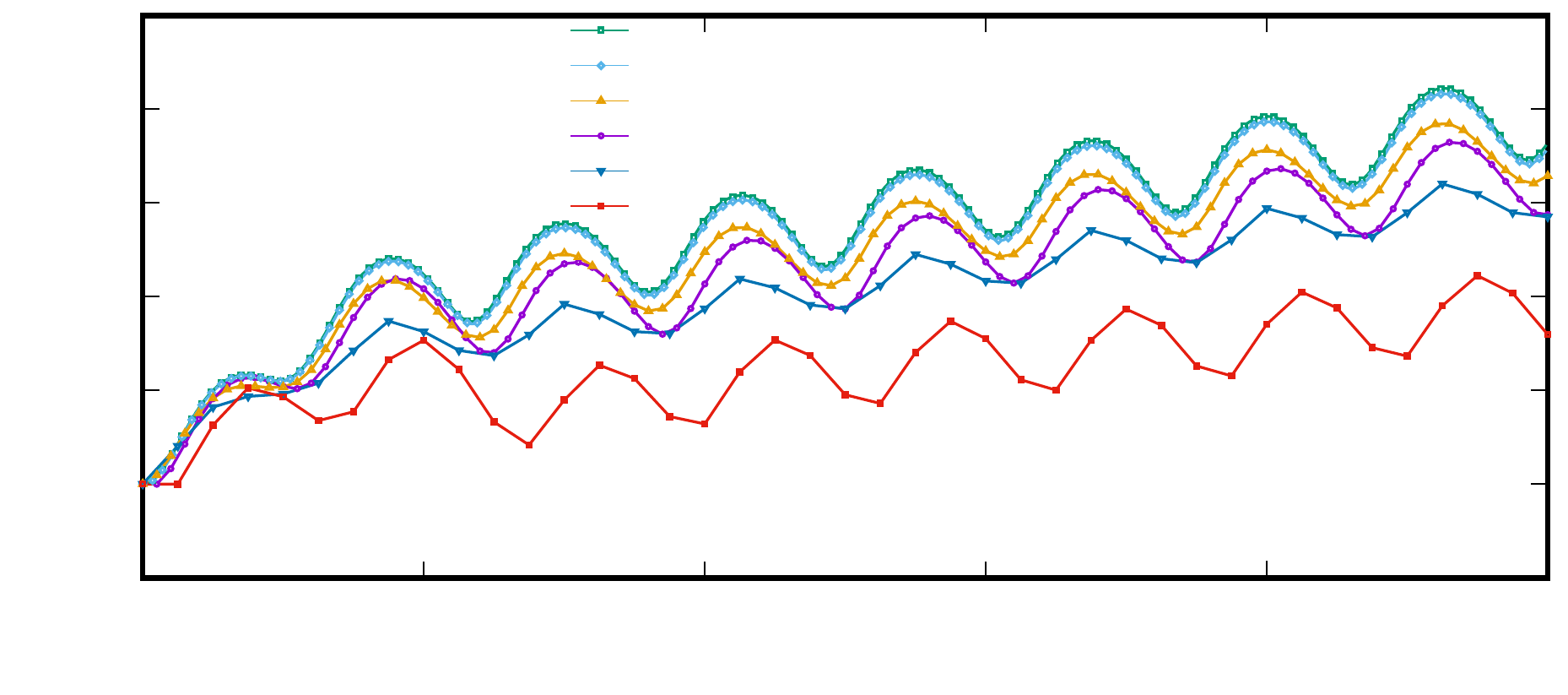}
       \caption{Swimmer position as a function of time for 
       the explicit, pseudo-implicit and semi-implicit and different time steps.}
       \label{fig:compa_x1c_methods}
   \end{center}
\end{figure}

A sensitive variable to compare is the viscous dissipation $\int_{\Omega_f} 2\mu\,\nabla^S{\mathbf{u}_h} : \nabla^S{\mathbf{u}_h}$, which is shown in
Fig. \ref{fig:compa_dissip_methods} 
as a function of time for all methods and different time steps. The explicit method is unstable for $\delta{t} > 10^{-3}$,
while in the others we can take $\delta{t}$ much larger without algorithmic crash. All methods converge to the same solution as $\delta t \to 0$. The semi-implicit and pseudo-implicit methods provide good approximations for $\delta t = 10^{-2}$. However, notice that the pseudo-implicit method requires several non-linear iterations
(tipically 3 to 5) to achieve convergence, which the semi-implicit, being a linear scheme, does not require. Results are quite poor for both,
the semi-implicit and pseudo-implicit methods for $\delta{t} = 2.5 \times 10^{-2}$.
The comparison is complemented with Fig. \ref{fig:compa_x1c_methods}, 
which assesses the impact of the temporal scheme and time step on the 
swimmer's advance. The explicit and semi-implicit
schemes for $\delta{t} = 10^{-3}$ deliver very similar results (green squares and
light blue diamonds, respectively). Finer time steps produce results that are undistinguishable
from those. The semi-implicit and the pseudo-implicit schemes for $\delta{t} = 10^{-2}$ (yellow triangles 
and magenta circles, respectively) produce reasonably accurate results, although some differences can be noticed. There appears a phase error in the advance of the swimmer while the average velocity after the initial transient is essentially correct. 
For $\delta{t} = 2.5\times 10^{-2}$ the error worsens, 
but still the main features of the swimmer's movement persist, however, 
the semi-implicit method (red solid squares) clearly exhibits a more dissipative behavior 
than the pseudo-implicit one (blue solid triangles). 

Summarizing, the semi-implicit method with intermediate time steps is a suitable and non-expensive choice 
that captures the swimming of finite length flexible rods in viscous fluids.

\begin{figure} 
   \begin{center}
       \scalebox{0.095}{\input{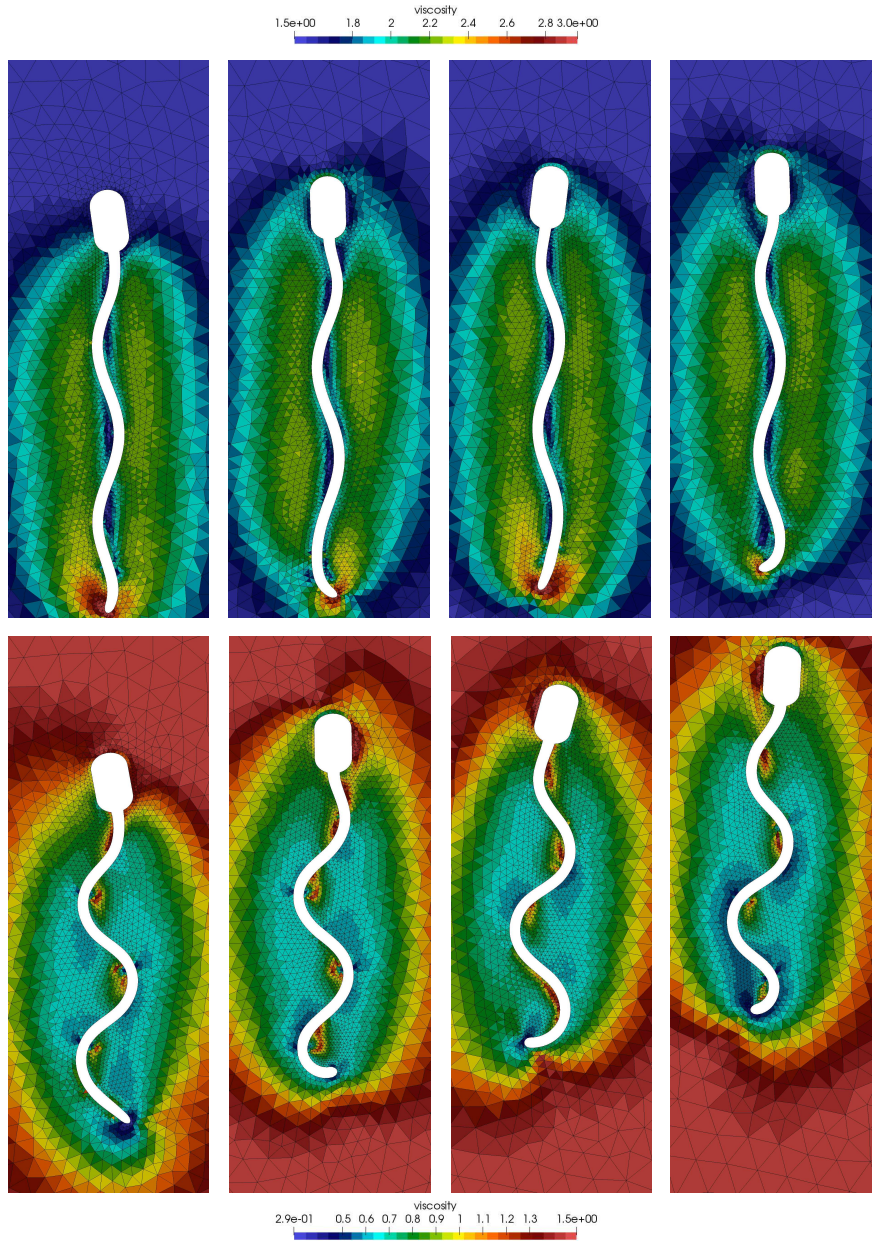}}
       \caption{Contours of viscosity and mesh for a swimmer 
       in a Carreau-Yasuda fluid at different times (from left to right,
       $t = 1.0,~1.8,~2.4,~3.2$). The top row corresponds to shear thickening 
       (power index $r = 1.15$) and the bottom row to shear thinning 
       (power index $r = 0.7$). The rod is discretized with $H = 1/8$ and 
       the fluid with typical refinement $h = 1/100$.}\label{fig:swimmer_carreau}
   \end{center}
\end{figure}

To conclude, we consider the case of a non-Newtonian fluid behavior. 
To that end, a Carreau-Yasuda \cite{MontenegroJohnson2012, Kennedy2013} viscosity law is considered,
namely,
\begin{equation}
 \mu(\dot{\gamma}) = \eta_{\infty} + (\eta_0 -  \eta_{\infty} ) (1 + 2\,\lambda \, \dot{\gamma}^2)^{\frac{r-1}{2}}
\end{equation}
where $\dot{\gamma} = \sqrt{\nabla^S{\mathbf{u}} : \nabla^S{\mathbf{u}}}$, 
$\eta_0 = 1.5$, $\eta_{\infty} = 10^{-3}$, $\lambda = 1$ and $r$ is a power
index. A shear thinning or shear thickening behavior is acquired depending on
whether $r < 1$ or $r > 1$, respectively. The swimmer is shown at different times
in Fig. \ref{fig:swimmer_carreau}. The bottom part corresponds
to the shear thinning case with $r = 0.7$, whilst the top part corresponds to the shear thickening case
with $r = 1.15$. The rod is discretized with $8$ elements and the fluid mesh with elements of 
characteristic size $h=1/100$ close to the wet boundary. 
Fig. \ref{fig:swimmer_carreau} also shows the fluid 
viscosity, which is an elementwise constant field, since linear elements are being used
for discretization of the Stokes problem. 
Note that the swimmer exhibits a significantly 
larger deformation and bigger net displacement in the shear thinning fluid as compared 
to the shear thickening one. Nevertheless, the aim at this point is not to perform a fair comparison 
of the swimmers' performance attained in each case, since the effective viscosity 
can be quite different, but to show the complex viscosity patterns that appear 
as the swimmers deform, which are plainly handled by the proposed formulation.
At this point, let us recall that the fluid nonlinearity is solved at each 
evaluation of the residual function by means of a Newton-Raphson procedure
with line search available in the SNES PETSc through the Firedrake interface.

%%% Local Variables:
%%% mode: latex
%%% TeX-master: "main"
%%% End:

\section{Concluding remarks}

This work proposes a fluid-structure interaction framework for Cosserat rods immersed in complex flows. The inertia of the rod and of the fluid are assumed negligible, as is typical in microsystems. Such an approach is suitable to study soft-bio-matter problems, e.g., elongated microorganism or appendices such as flagella and/or cilia. To showcase the promising features of the proposed approach, we specialize the new framework to the case of a planar non-shearable rod surrounded by an incompressible fluid with generalized Newtonian rheology.

As the chosen mechanical model relies on a configuration space that avoids shearing by construction, the deformation measures as well as the position of the wet surface depend only on the centroidal curve of the rod. Due to that, a discrete version of the model (in the sense of the finite element method) can be easily built and handled within the Firedrake platform, taking full advantage of its capabilities for modeling complex systems. The source code, which is made available alongside this manuscript, is lean and quite intuitive. The simplicity of the proposed framework for modeling active distributed devices/organisms, as needed to investigate locomotion, are noteworthy. 

The flexibility of Firedrake for modeling complex fluids allows to investigate fluid-structure interactions involving shear-thinning or shear-thickening Newtonian models. Nonlinear constitutive behaviors are rarely addressed in the literature and cannot be treated with other methods such as boundary elements. The fluid's rheology, however, is sometimes essential to model elongated microorganisms immersed in biological fluids \cite{RodrigoVelezCordero2012}. 

Though the reported work is rather extensive, there are many possible directions for further developments. The interaction among several elongated appendices as well as the optimization of swimming trajectories are just a few examples among a plethora of interesting applications that can be addressed in future works. The code is made freely available for interested colleagues to download and adapt to their own goals.

%%% Local Variables:
%%% mode: latex
%%% TeX-master: "main"
%%% End:

\section*{CRediT author statement}

\textbf{Roberto Federico Ausas:} Conceptualization, Methodology, Software, Investigation, Writing - Original Draft, Writing - Review \& Editing.
\textbf{Cristian Guillermo Gebhardt:} Conceptualization, Methodology, Formal analysis, Investigation, Writing - Original Draft, Writing - Review \& Editing.
\textbf{Gustavo Carlos Buscaglia:} Conceptualization, Methodology, Formal analysis, Investigation, Writing - Original Draft, Writing - Review \& Editing.

\section*{Declaration of Competing Interest}

The authors declare that they have no known competing financial interests or personal relationships that could have appeared to influence the work reported in this paper.

\section*{Acknowledgement}

Roberto Federico Ausas and Gustavo Carlos Buscaglia gratefully acknowledge financial support from the São Paulo Research Foundation FAPESP and from the Conselho Nacional de Desenvolvimento Científico e Tecnológico (grants CEPID-CeMEAI 2013/07375-0 and INCT-MACC).

Finally, the authors thank to D. Ham and K. Sagiyama from the Firedrake project and P. Farrell for giving some guidelines in the initial stages of the development.

\bibliographystyle{unsrt}

\end{document}